\pdfoutput=1
\documentclass[11pt]{article}

\usepackage[margin=1in]{geometry}
\usepackage[utf8]{inputenc}
\usepackage{url}
\usepackage{booktabs}
\usepackage{amsfonts}
\usepackage{graphicx}
\usepackage{amsmath,amsfonts,amssymb,amsthm}
\usepackage{subcaption}
\usepackage{enumitem}
\usepackage{setspace}
\usepackage{lmodern}
\usepackage{caption}
\usepackage{xspace}
\usepackage{hyperref} % Last

\theoremstyle{plain}
\newtheorem{theorem}{Theorem}[section]

\theoremstyle{definition}
\newtheorem{definition}[theorem]{Definition}
\newtheorem{assumption}[theorem]{Assumption}
\theoremstyle{remark}
\newtheorem{remark}[theorem]{Remark}

\newcommand{\Poincare}{Poincar\'e\xspace}

\def\keywordname{{\bfseries \emph{Keywords}}}%
\def\keywords#1{\par\addvspace\medskipamount{\rightskip=0pt plus1cm
\def\and{\ifhmode\unskip\nobreak\fi\ $\cdot$
}\noindent\keywordname\enspace\ignorespaces#1\par}}

\newtheoremstyle{plain}
  {0.5\baselineskip}
  {0.5\baselineskip}
  {\itshape}
  {}
  {\bfseries}
  {.}
  {0.5em}
  {}
\theoremstyle{plain}

\date{}
\title{Pancreatic $\beta-$Cell Dynamics with Three-Time-Scale Systems}
\author{
  \href{https://orcid.org/0000-0003-4207-3032}{Navojit Dhali Pallab} \\
  Mathematical Institute, Tohoku University, Sendai,980-8578, Japan\\
  \texttt{pallab.navojit.dhali.t4@dc.tohoku.ac.jp}
}

\begin{document}
\normalfont
\maketitle

\begin{abstract}
Pancreatic $\beta-$cells regulate insulin secretion through complex oscillations, which are vital for glucose control and diabetes research. In this paper, an existing mathematical model of $\beta-$cell dynamics is analyzed using a three-time-scale framework to study interactions among fast, intermediate, and slow variables. Through Geometric Singular Perturbation Theory (GSPT), the influence of ATP on oscillatory dynamics via membrane potential is explored. At the non-hyperbolic point, where standard methods fail, blow-up analysis is applied to investigate canard dynamics shaped by intermediate and slow variables. Numerical simulations with varied parameters reveal the glucose-dependent oscillations linked to slow dynamics near the pseudo-singular points. By leveraging the pseudo-singular point, the linger time is defined, and simulated results for the coupling strength needed for bursting initiation synchronization are presented as a sufficient condition. This study links mathematics and biology, offering insights into diabetic studies.
\end{abstract}

\def\keywordname{\textbf{Keywords:}}
\def\keywords#1{\par\addvspace{0.5\baselineskip}\noindent\keywordname\ \textit{#1}\par\vspace{0.5\baselineskip}}

\section{Introduction}
Insulin is the only glucose-lowering hormone in the body \cite{rorsman2003insulin}. Pancreatic $\beta-$cells regulate blood glucose by secreting insulin in accordance with the blood glucose concentration.  The electrical activity of the Pancreatic $\beta-$cells exhibits bursting behavior, which leads to oscillations in insulin secretion \cite{gilon1993oscillations} led by oscillations in the intercellular free $Ca^{2+}$ oscillation \cite{pernarowski1998fast, bertram2008phantom, bertram2010electrical}. Because of the central importance of glucose homeostasis, the study of the Pancreatic $\beta-$cell has had a remarkable theoretical interest for the last few decades, besides experimental study \cite{bertram2000phantom,pedersen2005intra}.

Complex systems in the real world can be modeled as a network of connected components, and synchronization is one of the fundamental phenomena in the study of network dynamics. For the study of the dynamics of the Pancreatic $\beta-$cells network, its importance lies in the key role of accelerating the insulin secretion because of coordinated bursting. This complex metabolic process in the Pancreatic $\beta-$cell critically depends on ATP-dependent potassium channels by modulating the membrane potential. In order to understand the importance of the K(ATP) driven oscillation of slow calcium oscillation, Marinelli et. al. \cite{marinelli2022oscillations} proposed a Mathematical model, which incorporate (i) membrane potential, (ii) activation variable for delayed rectifier of potassium ion channel, (iii) free cytosolic calcium concentration, (iv) free calcium concentration in the endoplasmic reticulum, and (v) ADP/ATP ratio as the state variables. In this study, we focused on the dynamics of the Pancreatic $\beta-$cells of this mathematical model from a viewpoint of three-time-scale dynamical systems with two fast, one intermediate, and two slowest variables. At first, we will present the network of the $N$ coupled Pancreatic $\beta-$cells for the mathematical model adopted from \cite{marinelli2022oscillations}.

Consider a network of $N-$number of heterogeneous Pancreatic $\beta-$cells that are electrically connected (or, connected via gap junction) through a coupling function in their fast state variable with coupling strength $k$, where $k \in \mathbb{R^+}$. Let $V_i=V_i (t) \ (mV)$ be the membrane potential of the $i^th$ Pancreatic $\beta-$cell at time $t\ (ms)$ (in this model, we use $ms$ as time unit for all the variables and subscript $i$ of state variables to represent the $i^{th}$ cell of the network), and $C_i=C_i (t)$  $(mM)$ and $E_i=E_i (t)$  $(mM)$ be the calcium ion concentration in cytosolic plasma membrane, endoplasmic reticulum of this cell, respectively. Consider $K_i=K_i (t)\in [0,1]$ as the voltage ($V_i (t)$) dependent potassium ion ($K^+$) gating variable and $R_i=R_i (t)$ as the ADP/ATP ratio in the intracellular $i^{th}$ Pancreatic $\beta-$cell.

There are many choices for how we could define the coupling function $H_i (V_1,\cdots, V_N;\epsilon):\ \mathbb{R}^N \to \mathbb{R}$, for instance, coupling between two cells can be represented as the difference function of their state variables that could be linear or nonlinear. 
With the above description of state variables and arbitrary coupling functions, we can write the canonical network model for the Pancreatic $\beta-$cell as
\begin{equation}
    \begin{aligned}
        \frac{dV_i}{dt} &= F_{1i} (V_i,K_i, C_i, E_i,R_i,\mu_i )+k H_i (V_1,V_2,\cdots,V_N;k),\\
        \frac{dK_i}{dt} &= \frac{1}{\tau_K} F_{2i} (V_i,K_i, C_i, E_i,R_i,\mu_i ),\\
        \frac{dC_i}{dt} &= F_{3i} (V_i,K_i, C_i, E_i,R_i,\mu_i),\\
        \frac{dE_i}{dt} &= F_{4i} (V_i,K_i, C_i, E_i,R_i,\mu_i),\\
        \frac{dR_i}{dt} &=\frac{1}{\tau_R}  F_{5i} (V_i,K_i, C_i, E_i,R_i,\mu_i),
    \end{aligned}\label{nbeta}
\end{equation}
for \(i=1,2,\cdots ,N\). The explicit expression of the functions used in (\ref{nbeta}) is as follows
\begin{equation*}
    \begin{aligned}
        F_{1i} (V_i,K_i, C_i, E_i,R_i,\mu_i) &=-\frac{1}{c_v}  [g_{Ca}  M_{\infty} (V_i )  (V_i-V_{Ca} )\\ &  +(g_K K_i+g_{KCa}  \frac{C_i^5}{C_i^5+K_D^5 }+g_{KATP}  R_i )(V_i-V_K )],\\
         F_{2i} (V_i,K_i, C_i, E_i,R_i,\mu_i ) &=K_{\infty} (V_i )-K_i,\\
        F_{3i} (V_i,K_i, C_i, E_i,R_i,\mu_i) &=- a_1 g_{Ca}  M_{\infty} (V_i )  (V_i-V_{Ca}) - a_2 C + a_3 E_i + a_4,\\
        F_{4i} (V_i,K_i, C_i, E_i,R_i,\mu_i ) &= a_6  C_i - a_5  (E_i-C_i),\\
        F_{5i} (V_i,K_i, C_i, E_i,R_i,\mu_i ) &= R_{\infty} (C_i )-R_i.
    \end{aligned}\label{nbetafunc}
\end{equation*}
Here, $\mu_i$ is the parameters set of the $i^{th}$ Pancreatic $\beta-$cell. By heterogeneous, we mean that the values of parameter sets are nonidentical. For simplicity of representation, we have not used the subscript for the parameters in the explicit expression of functions $F_i$. Besides these, we use the Boltzmann function (another interpretation of these functions: fast threshold modulation form \cite{belykh2005synchronization} are used as coupling function) to define the saturated value of voltage-dependent calcium and potassium gating variable, and ADP/ATP ratio with$ M_{\infty} (V_i ), K_{\infty} (V_i )$ and $R_{\infty} (C_i)$, respectively, as
\begin{equation}
    \begin{aligned}
        M_\infty (V_i) = \frac{1}{1+\exp{(\frac{v_M-V_i}{s_M})}},\\
        K_\infty (V_i) = \frac{1}{1+\exp{(\frac{v_K-V_i}{s_K})}},\\
        R_\infty (C_i) = \frac{1}{1+\exp{(\frac{v_R (G)-C_i}{s_R})}},
    \end{aligned}\label{nbetaBF}
\end{equation}
where $v_R (G):=\frac{G-p_R}{k_R} $ is the function of plasma glucose concentration, $G\ (mM)$. The physiological description of the parameters used in the model and Boltzmann functions is presented in Table \ref{tebParaOG}. 

Since the network of N coupled $\beta-$cells is the system of $5N$ coupled ordinary differential equation in (\ref{nbeta}), synchronization between these cells in the network is relevant to finding a stable limit cycle for the fast variables $V_1, V_2,\cdots, V_N$ for some $k_0>0, k\in (0,k_o ]$  for $i=1,2,\cdots, N$.

The remaining part of this paper is organized as follows: Section~\ref{sec:singleBETAcell} presents the three-time-scale systems derived from the model with the critical manifold analysis, Section~\ref{sec:canard} details findings on the fold singularity, normal form of the model at that point and its dynamics, Section~\ref{sec:numerical} discusses the limit cycle and relevant dynamics of the systems, and Section~\ref{sec:conclusion} presents the conclusion of this study. 

Before continuing to the next section, we want to define a few well-known definitions here.
\begin{definition}[Limit Cycle]\label{def:limitcycle}
    Consider a two-dimensional dynamical system of the form \[ \dot{x}(t) = f(x(t)), \] where $f:\mathbb{R}^2 \rightarrow \mathbb{R}^2$ is a smooth function. A trajectory of this system is some smooth function $x(t):=[x_1(t),x_2(t)]^T$ with values in $\mathbb{R}^{2}$ which satisfies this differential equation. Such a trajectory is called closed (or periodic) if it is not constant but returns to its starting point, i.e. if there exists some $T_0>0$ such that  $x(t+T_0)=x(t)$ for all $t\in \mathbb{R}$. A limit cycle is a closed trajectory which is the limit set of some other trajectory. Define $(x_1(t),x_2(t))=\gamma(t)$. The limit cycle $\Gamma$ can be expressed as \[\Gamma:=\{\gamma(t)\ |\ 0\leq t <T_0\},\] where $T_0$ is the minimal period of limit cycle.
\end{definition}

\begin{definition}[\Poincare Map \cite{wiggins2003intro}]
Consider a dynamical system
\begin{equation} \label{sysX}
    \begin{aligned}
        \dot{X}(t)=F(X),\ \ \ \ \ X\in\mathbb{R}^n,
    \end{aligned}
\end{equation}
where $F:U\rightarrow\mathbb{R}^n$ is $C^r$ on some open set $U\subset\mathbb{R}^n$. Let $\phi(t,\cdot)$ denote the flow generated by (\ref{sysX}). Suppose that (\ref{sysX}) has a periodic solution of period $T$ which we denote by $\phi(t,X_0)$, where $X_0 \in \mathbb{R}^n$ is any point through which this periodic solution passes.\\
\\Let $S$ be an $(n-1)$ dimensional surface transverse to the vector field at $X_0$. In particular, $S$ is a cross-section to the vector field (\ref{sysX}), called the \textit{\Poincare section}. Since $\phi(t,X)$ is $C^r$ (because of $F(X)$ is $C^r$), one can find an open set $V\subset S$ such that the trajectories starting in $V$ return to $S$ in a time close to $T$. The map that associates points in $V$ with their points of first return to $S$ is called the \textit{\Poincare map}, which denoted by $P$,
\begin{equation*}
    \begin{aligned}
        P:V\rightarrow S,\\
        X\mapsto\phi(\tau(X),X),
    \end{aligned}
\end{equation*}
where $\tau(X)$ is the time of first return of the point $X$ to $S$. Here, by construction, $\tau(X_0)=T$ and $P(X_0)=X_0$.\\
\\A fixed point of $P$ corresponds to a \textit{periodic orbit} of (\ref{sysX}), and a \textit{period $k$ point} of $P$ corresponds to a periodic orbit of (\ref{sysX}) that pierces $S$ for $k$ times before closing,$i.e.$, for a point $X\in V$ such that $P^k(X)=X$ provides $P^j(X)\in V,\ j=1,2,\cdots,k$.
\end{definition}

\begin{definition}[Canard solution \cite{tchizawa2007generic}]\label{def:canard}
 Consider a system 
 \begin{equation}\label{SForm}
    \begin{aligned}
        \frac{dx}{dt} &= \tilde{f}(x,y,\varepsilon),\\
        \frac{dy}{dt} &= \varepsilon \tilde{g}(x,y,\varepsilon),       
    \end{aligned}
\end{equation}
where $x\in \mathbb{R}$ and $y\in \mathbb{R}$ are the fast and slow variables with $1\leq k < n$, respectively, and the vector valued functions $\Tilde{f}$ and $\Tilde{g}$ are sufficiently smooth. A solution $(x(t,\varepsilon), y(t,\varepsilon))$ of the systems are called canard, if there exist $t_1 < t_0 < t_2$ such that $(i)$ standard part of $(x(t_0, \varepsilon), y(t_0,\varepsilon))$ belongs to the critical manifold $\mathcal{C}=\{ (x,y)\in\mathbb{R}^2 \mid \ f(x,y,0)=0 \}$, $(ii)$ for $t\in (t_1,t_0)$ the segment of the trajectory $(x(t,\varepsilon), y(t,\varepsilon))$ is infinitesimally close to the attracting branch of $\mathcal{C}$, $(iii)$ for $t\in (t_0, t_2)$ the segment of trajectory remain close to the repelling branch of $\mathcal{C}$, and $(iv)$ the attracting and repelling parts of the trajectory are not infinitesimally small. 
More concisely, a trajectory segment of a fast-slow system (\ref{SForm}) is a canard if it stays within $\mathcal{O}(\varepsilon)$ distance to a repelling branch of a slow manifold for a time that is $\mathcal{O}(1)$ on the slow time scale $\tau= \varepsilon t$.
\end{definition}

\begin{table}[ht]
\centering
    \begin{tabular}{|p{2cm}|p{11.75cm}|}
        
        \hline
        Parameter & Description \\
        \hline
        $c_v$ & membrane capacitance\\
        $g_{Ca}$ & maximal $Ca^{2+}$ conductance through voltage-dependent $Ca^{2+}$ channels\\
        $V_{Ca}$ & equilibrium (Nernst) membrane potential for $Ca^{2+}$\\
        $g_K$ & maximal $K^+$  conductance through voltage-gated $K^+$ channels\\
        $g_{KCa}$ & maximal $K^+$  conductance through $Ca^{2+}$-dependent $K^+$ channels\\
        $K_D$ & equilibrium dissociation constant\\
        $g_{KATP}$ & maximal $K^+$  conductance through ATP-dependent $K^+$ channels\\
        $V_K$ & equilibrium (Nernst) membrane potential for $K^+$\\
        $a_1$ &	the ratio between $Ca^{2+}$  concentration and $Ca^{2+}$ current\\
        $a_2,a_3, a_5, a_6$ &	$Ca^{2+}$ exchange rate among cytosolic and endoplasmic compartment\\
        $a_4$ &	difference of constantly cytosolic $Ca^{2+}$ leakage and $Ca^{2+}$ exchange\\
        $\tau_K$ &	time constant for $K^+$ gating variable\\
        $\tau_R$ &	time constant for ADP/ATP ratio \\
        $v_M$ &	the membrane potential at which $Ca^{2+}$ gating probability reaches its half of the maximum opening ($v_M$ is the inflection point of the probability curve)\\
        $s_M$ &	slope of steady-state activation for $Ca^{2+}$ gating variable\\
        $v_K$ &	the membrane potential at which $K^+$ gating probability reaches its half of the maximum opening\\  
        $s_K$ & slope of steady-state activation for $K^+$ gating variable\\
        $s_R$ &	slope of steady-state activation for ADP/ATP ratio\\
        $p_R, k_R$ &	$v_R (G)$  is the cytosolic $Ca^{2+}$ concentration, $C$, when the ADP/ATP ratio is at its half-maximum value, $i.e.$, parameter $p_R$ and $k_R$ are glucose threshold and glucose sensitivity scaling of the $\beta-$cell, respectively\\
        \hline
    \end{tabular}
    \caption{Physiological description of parameters used for network model (\ref{nbeta}). For more details of the parameters, see \cite{marinelli2022oscillations, pedersen2009newcomer, wang2008identifying, bertram2004filtering}.}
    \label{tebParaOG}
\end{table}

\begin{table}[ht]
\centering
    \begin{tabular}{|p{1.75cm}|p{4.75cm}||p{1.75cm}|p{4.75cm}|}
        
        \hline
        Parameter & Value & Parameter & Value  \\
        \hline
        $c_v$ & $5300$ $fF$ & $G$ & $8$ $mM$\\
        $\tau_K$ & $16$ $ms^{-1}$ &$\tau_R$ & $3 \times 10^5$  $ms^{-1}$\\
        $g_{Ca}$ & $1200$ $pS$ & $V_{Ca}$ & $25$ $mV$\\
        $g_K$ & $3000$ $pS$ & $g_{KCa}$ &$10$ $pS$\\
        $g_{KATP}$ & $500$ $pS$ & $K_D$ & $0.3$ $\mu M$\\
        $V_K$ & $-16$ $mV$ & $f_C$ &	$0.01$\\
        $f_{er}$ & $0.01$ & $v_{rat}$ &  $5$ \\
        $\alpha$ & $4.5 \times 10^{-6}$ $fA^{-1} \mu M\ ms^{-1}$	& $k_{PMCA}$ & $0.2$ $ms^{-1}$	\\
        $k_{SERCA}$ &	$0.4$ $ms^{-1}$ & $p_{leak}$ & $0.0005$ $ms^{-1}$	\\
        $s_M$ &	$12$ $mV$ & $v_M$ & $-20$ $mV$\\  
        $s_K$ & $5$ $mV$ &  $v_K$ & $-16$ $mV$\\
        $s_R$ & $0.1$	$\mu M$& $p_R$ & $1.75$ $mM$	\\
        $k_R$ & $58$ $mM$ & & \\
        \hline
    \end{tabular}
    \caption{Parameter values used for (\ref{ScaledMarinelli}). For more details about the parameters, see \cite{marinelli2022oscillations, wang2008identifying, bertram2004filtering}.}
    \label{tebPara}
\end{table}

\begin{figure}[ht]
    \centering
    \includegraphics[width=0.6\textwidth]{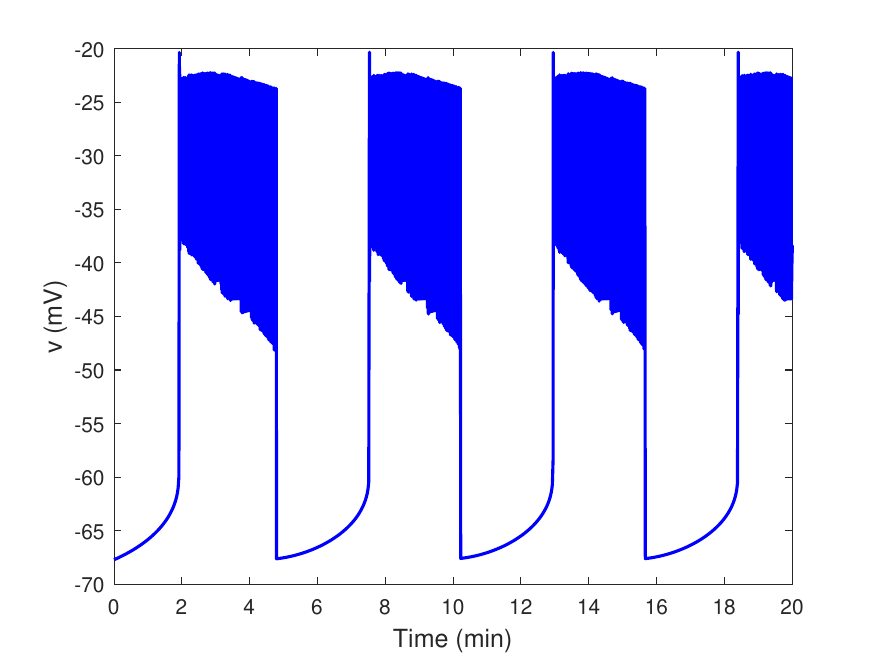}
    \caption{Bursting oscillation of simplified $\beta-$cell model (\ref{ScaledMarinelli}) of Marinelli et al. model \cite{marinelli2022oscillations}, for $G=13\ mM$.}
    \label{fig:single_beta_burst}
\end{figure}

\section{Three-Time-Scale Pancreatic Beta-Cell Model}\label{sec:singleBETAcell}
Consider the Pancreatic $\beta-$cell activity model corresponds to the K(ATP) conductance drive slow calcium oscillations \cite{marinelli2022oscillations} as the following form
\begin{equation}\label{ScaledMarinelli}
\begin{aligned}
    \frac{dv}{dt} &= h_1 (v,u,x,y),\\
    \frac{du}{dt} &= h_2 (v, u),\\
    \frac{dx}{dt} &= \varepsilon f (v, x, z),\\
    \frac{dy}{dt} &= \varepsilon \delta  g_1 ( x, y),\\
    \frac{dz}{dt} &= \varepsilon \delta  g_2(x, z),
\end{aligned}
\end{equation}
where $v,u,x,y,z \in \mathbb{R}$, are equivalent to $V_i, K_i, C_i, E_i, R_i$, respectively, for $N=1$, and, $h_1, h_2, f, g_1$ and $g_2$ are defined as
\begin{equation}
    \begin{aligned}
        h_1 (v,u,x,y) &:= -\left(a_1 X_\infty (v) (v-a_2) + \left(a_3 u + a_4 U_\infty (x) + a_5 y \right)(v-a_6) \right),\\
        h_2 (v,u) &:= a_7 \left(Y_\infty (v) -u \right),\\
        f (v, x, z) &:= -d_1 X_\infty (v) (v-a_2) - d_2 x + d_3 z,\\
        g_1(x,y;G) &:= Z_\infty (x) -y ,\\
        g_2 ( x, z) &:= k_1 x - k_2 z,\\
    \end{aligned}
\end{equation}
with
\begin{equation*}
    \begin{aligned}
        X_\infty (v) &= \left[1+exp\left(\frac{v_1-v}{s_1}\right) \right]^{-1},\\
        Y_\infty (v) &= \left[1+exp\left(\frac{v_2-v}{s_2}\right) \right]^{-1},\\
        Z_\infty (x) &= \left[1+exp\left(\frac{\frac{G-p_r}{k_r}- x}{s_3}\right) \right]^{-1},\\
        U_\infty (x) &= \left[1+\left(\frac{K_d}{x}\right)^5\right]^{-1},\\
    \end{aligned}
\end{equation*}
In system (\ref{ScaledMarinelli}), $v$ and $u$ are the fast variables, $x$ is the intermediate variable, $y$ and $z$ are the slowest variables. All the parameter values for this system are presented in (Table~\ref{tab:ScalePara}), which are calculated using the parameter values mentioned in \cite{marinelli2022oscillations}.

\begin{figure}[ht]
    \centering
    \begin{subfigure}[t]{0.49\textwidth}
        \centering
        \begin{subfigure}[t]{\textwidth}
           \includegraphics[width=\linewidth, height=5.5cm]{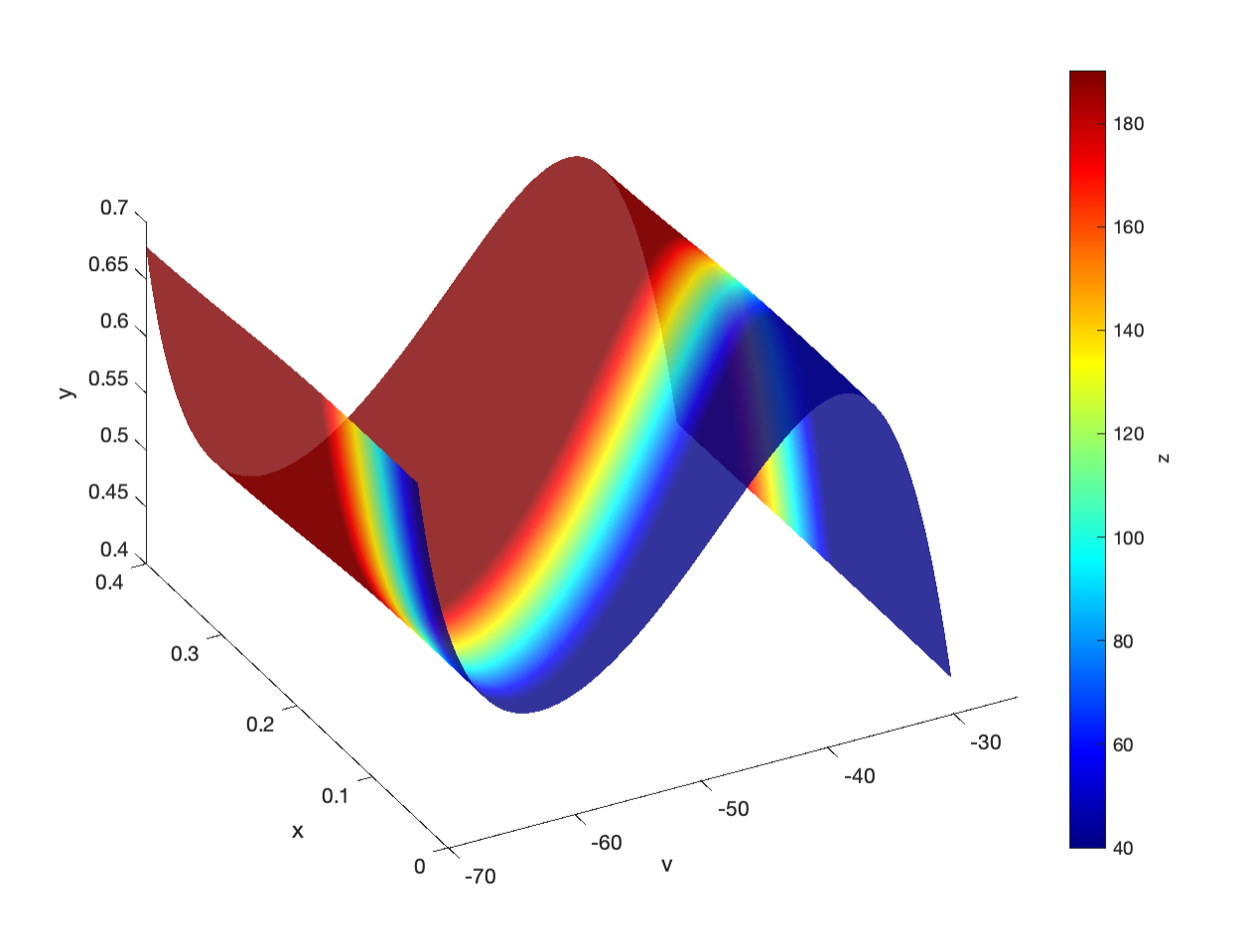}
           \caption{Critical manifold for $h_1 = h_2 =0$.}
           \label{fig:h1h20}
        \end{subfigure}
    \end{subfigure}
    \hfill 
    \begin{subfigure}[t]{0.49\textwidth}
        \centering
        \begin{subfigure}[t]{\textwidth}
           \includegraphics[width=\linewidth, height=5.5cm]{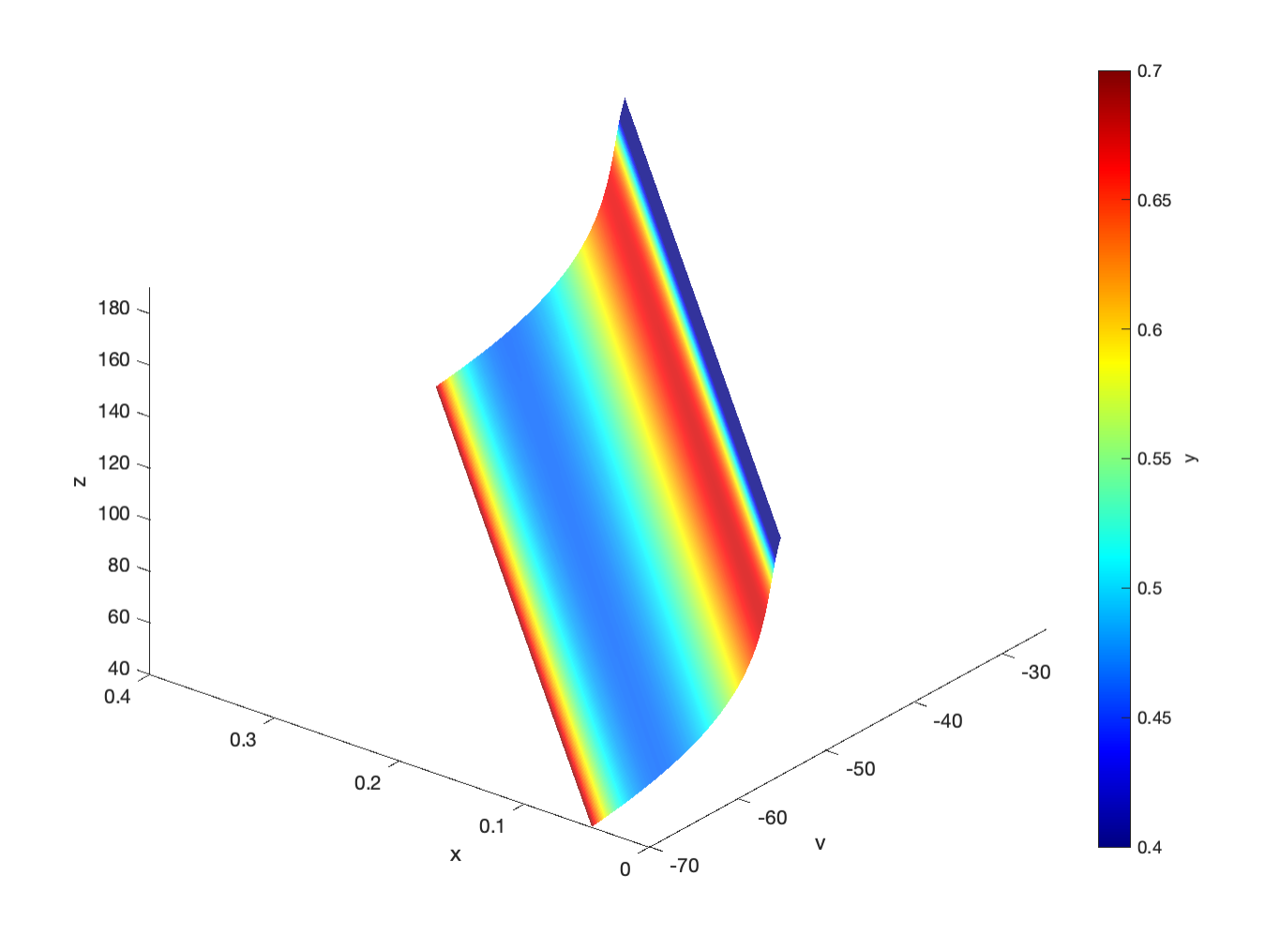}
           \caption{Critical manifold for $h_1=h_2=f=0$ in addition to (a).}
           \label{fig:h1h2f0}
        \end{subfigure}
    \end{subfigure}
    \caption{Critical manifold of the system (\ref{ScaledMarinelli}), (a) in ($v,x,y$) and (b) in ($v,x,z$). The color bar for $z$ in (a) shows how the critical manifold $\mathcal{C}_1$ (\ref{cm:eps0delta0}) embedded on $\mathcal{C}$ (\ref{cm:eps0}).}
    \label{fig:criticalmani}
\end{figure}

\begin{figure}[ht]
    \centering
    \includegraphics[width=  0.7\textwidth, height=13cm, angle=0]{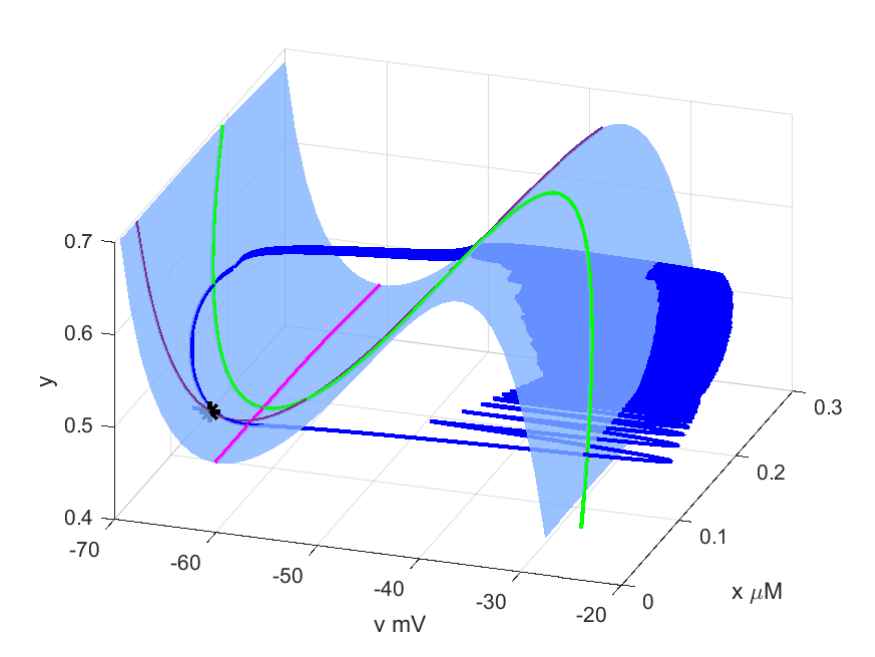} 
    \caption{Critical manifold ($\varepsilon$=0 and $\delta=0$) and the trajectory ($t\in[0,30]$ minutes) of the full system (blue) for $G=8\ mM$, and other parameters remain the same as Table \ref{tab:ScalePara}. When the trajectory reaches the fold curve $\mathcal{L}$ (\ref{cm:eps0foldcurve}) (magenta) of the critical manifold $\mathcal{C}$, it leaves the stable branch.}
    \label{FIG:TrajectoryWITHcriticalmanifold}
\end{figure}

\begin{table}[ht]
    \centering
    \begin{tabular}{|c|c||c|c|}
    \hline
      $a_1$   & 0.226 $ms^{-1}$& $a_2$ & 25 $mV$\\
      $a_3$& 0.566038 $ms^{-1}$& $a_4$  & 0.00189 $ms^{-1}$ \\
      $a_5$  & 0.0943 $ms^{-1}$& $a_6$ & -75 $mV$\\
      $a_7$ & 0.0625 $ms^{-1}$& $v_1$ & -20 $mV$ \\
      $v_2$ & -16 $mV$& $s_1$ & 12 $mV$\\
      $s_2$ & 5 $mV$& $s_3$ & 0.1 $\mu M$\\
      $d_1$ & 0.144  $ms^{-1}$& $d_2$ & 16 $ms^{-1}$\\
      $d_3$ & 0.013 $ms^{-1}$& $k_1$ & 5977.6 $ms^{-1}$ \\
      $k_2$ & 7.463 $ms^{-1}$& $p_r$ & 1.75 $mM$\\
      $k_r$ & 58 $mM$& $K_d$ & 0.3 $\mu M$\\
      $\varepsilon$ & $3.753 \times 10^{-4}$  & $\delta$ & 0.0089\\
      \hline
    \end{tabular}
    \caption{Parameters value used in this work for the model (\ref{ScaledMarinelli}) are adapted from \cite{marinelli2022oscillations}. We will treat $G$ as the control parameter within a certain range corresponding to the Pancreatic $\beta-$cell physiology.}
    \label{tab:ScalePara}
\end{table}

The dynamics of the system (\ref{ScaledMarinelli}) are governed in fast time $t$. By rescaling time as $\tau=\varepsilon t$, the system can be written as 
\begin{equation}
\begin{aligned}
   \varepsilon \frac{dv}{d\tau} &= h_1 (v,u,x,y),\\
    \varepsilon \frac{du}{d\tau} &= h_2 (v, u),\\
    \frac{dx}{d\tau} &=  f (v, x, z),\\
    \frac{dy}{d\tau} &= \delta  g_1 ( x, y),\\
    \frac{dz}{d\tau} &= \delta  g_2(x, z),
\end{aligned}\label{sys:intermediate}
\end{equation}
which is called the intermediate system. For $\varepsilon=0$, the dynamics of the system (\ref{sys:intermediate}) governed by the intermediate problem
\begin{equation}
\begin{aligned}
    \frac{dx}{d\tau} &=  f (v, x, z),\\
    \frac{dy}{d\tau} &= \delta  g_1 ( x, y),\\
    \frac{dz}{d\tau} &= \delta  g_2(x, z),
\end{aligned}\label{sys:intermediateProb}
\end{equation}
on the $3-$dimensional set of equilibrium $\{(v,u,x,y,z)\in\mathbb{R}^5 \mid h_1=h_2=0\}$, is the critical manifold of the system (\ref{ScaledMarinelli}) (Figure~\ref{fig:criticalmani}(A)), and can be written as
\begin{equation}
    \begin{aligned}\label{cm:eps0}
        \mathcal{C} = \left\{(v,u,x,y,z) \in \mathbb{R}^5 \ \lvert \ a_1 X_\infty (v) (v-a_2) + \left(a_3 Y_\infty (v)  + a_4 U_\infty (x) + a_5 y \right)(v-a_6) = 0,\ u= Y_\infty(v) \right\},
    \end{aligned}
\end{equation}
and the fold curves on $\mathcal{C}$ where determinant of $D_{(v,u)}(h_1,h_2))$ is zero, takes the form
\begin{equation}
    \begin{aligned}\label{cm:eps0foldcurve}
        \mathcal{L} = \left\{(v,u,x,y,z) \in \mathcal{C} \ \lvert \ y = \frac{a_1}{a_5} \left[\frac{V_1 (v) X_\infty (v)}{s_1}  + 1 \right] X_\infty (v) + \frac{a_3}{a_5} \left[\frac{V_2 (v) Y_\infty (v) (v-a_6)}{s_2}  + 1 \right] Y_\infty (v) + \frac{a_4}{a_5} U_\infty (x) \right\}.
    \end{aligned}
\end{equation}
Again, by scaling the time as $\tau_s=\delta \tau = \varepsilon \delta t$, the system (\ref{ScaledMarinelli}) can be written as the slowest system 
\begin{equation}
\begin{aligned}
   \varepsilon \delta \frac{dv}{d\tau_s} &= h_1 (v,u,x,y),\\
    \varepsilon \delta \frac{du}{d\tau_s} &= h_2 (v, u),\\
     \delta \frac{dx}{d\tau} &=  f (v, x, z),\\
    \frac{dy}{d\tau} &=  g_1 ( x, y),\\
    \frac{dz}{d\tau} &=  g_2(x, z),
\end{aligned}\label{sys:slowest}
\end{equation}
which is called the slowest system. For $\delta=0, \ \varepsilon=0$, the dynamics of the system is determined by the reduced problem
\begin{equation}
\begin{aligned}
    \frac{dy}{d\tau_s} &=  g_1 ( x, y),\\
    \frac{dz}{d\tau_s} &=  g_2(x, z),
\end{aligned}\label{sys:slowestProblem}
\end{equation}
on the $2-$dimensional critical manifold $\{f=0\} \cap\mathcal{C}$, defined by
\begin{equation}\label{cm:eps0delta0}
    \begin{aligned}
    \mathcal{C}_1 = \left\{ (v,u,x,y,z) \in \mathcal{C} \ \lvert \ x= \frac{d_3 z - d_1 X_\infty (v) (v-a_2)}{d_2} \right\}.
    \end{aligned}
\end{equation}

\begin{remark}
    The systems (\ref{ScaledMarinelli}), (\ref{sys:intermediate}) and \ref{sys:slowest} are equivalent. Though we have already found three systems, there is another representation of the system (\ref{ScaledMarinelli}), by defining the time scale $\tau_{ss}= \delta t$. In that case, for $\delta=0$, the dynamics on critical manifold $\mathcal{C}_1$ (\ref{cm:eps0delta0}) is governed by 
    \begin{equation}
    \begin{aligned}
        \frac{dy}{d\tau_{ss}} &= \varepsilon  g_1 ( x, y),\\
        \frac{dz}{d\tau_{ss}} &= \varepsilon g_2(x, z).
    \end{aligned}\label{sys:SSlowestProblem}
    \end{equation}
    One of the main differences between systems (\ref{sys:slowestProblem}) and (\ref{sys:SSlowestProblem}) is that former separates the fast and the slowest variables, where the unperturbed syatem of (\ref{ScaledMarinelli}) is defined by 
    \begin{equation}\label{unperturbed2fast}
    \begin{aligned}
        \frac{dv}{dt} &= h_1 (v,u,x,y),\\
        \frac{du}{dt} &= h_2 (v, u),\\
        \frac{dx}{dt} &=0,\\
        \frac{dy}{dt} &= 0,\\
        \frac{dz}{dt} &= 0,
    \end{aligned}
    \end{equation}
and the later one separates the slowest variables and fast-intermediates variables, in that case the unperturbed system of (\ref{ScaledMarinelli}) becomes
    \begin{equation}\label{unperturbed3fast}
    \begin{aligned}
        \frac{dv}{dt} &= h_1 (v,u,x,y),\\
        \frac{du}{dt} &= h_2 (v, u),\\
        \frac{dx}{dt} &=\varepsilon g(v,x,z),\\
        \frac{dy}{dt} &= 0,\\
        \frac{dz}{dt} &= 0,
    \end{aligned}
    \end{equation}
    which is again a singular perturbed problem.
\end{remark}

\section{Dynamics Near the Pseudo Singular Point}\label{sec:canard}
For the reduced problem (\ref{sys:slowestProblem}), the normalized slow dynamics on $\mathcal{C}_1$ have the following form
\begin{equation}
    \begin{aligned}
        \frac{dv}{d_{\tau_n}} &= (\partial_u h_2 \cdot \partial_x f \cdot \partial_y h_1) g_1 - ( \partial_x h_1 \cdot \partial_u h_2 \cdot \partial_z f) g_2, \\
        \frac{du}{d_{\tau_n}} &= -(\partial_v h_2 \cdot \partial_x f \cdot \partial_y h_1) g_1 + ( \partial_x h_1 \cdot \partial_v h_2 \cdot \partial_z f) g_2,\\
        \frac{dx}{d_{\tau_n}} &= -( \partial_u h_2 \cdot \partial_v f \cdot \partial_y h_1) g_1 + ( ( \partial_v h_1 \cdot \partial_u h_2) - ( \partial_u h_1 \cdot \partial_v h_2)) \partial_z f \cdot g_2,\\
        \frac{dy}{d_{\tau_n}} &= -\det[D_{(v,u,x)}(h_1,h_2,f)] g_1,\\
        \frac{dz}{d_{\tau_n}} &= -\det[D_{(v,u,x)}(h_1,h_2,f)] g_2,
    \end{aligned}\label{normalized5Dsystem}
\end{equation}
where $\tau_n=-\det[D_{(v,u,x)}(h_1,h_2,f)] t_s$. The equations for fast variables $v$ and $u$ in the constrained system (system \ref{normalized5Dsystem} in time scale $\tau_s$) are derived according to the Fenichel theory \cite{fenichel1979geometric} that the critical manifold $C_1$ can be considered as locally invariant under the flow of the original system (\ref{ScaledMarinelli}). The negative sign used for the time scaling in the system (\ref{normalized5Dsystem}) ensures the same flow direction of the constrained system and the normalized intermediate problem.

%%%%%%%%%%%%%%%%%%%%%%%%%%%%%
\begin{definition}[Pseudo-singular point]
    Let $P=(v_p, u_p, x_p, y_p, z_p)\in \mathcal{C}_1$ be the point which satisfies the following condition
   \begin{equation}
    \begin{aligned}
        \left( \partial_u h_2 \cdot \partial_x f \cdot \partial_y h_1 \right) g_1 - \left( \partial_x h_1 \cdot \partial_u h_2 \cdot \partial_z f \right) g_2 =0,\\
        -\left( \partial_v h_2 \cdot \partial_x f \cdot \partial_y h_1 \right) g_1 + \left( \partial_x h_1 \cdot \partial_v h_2 \cdot \partial_z f \right) g_2 =0,\\
        -\left( \partial_u h_2 \cdot \partial_v f \cdot \partial_y h_1 \right) g_1 + \left[ \left( \partial_v h_1 \cdot \partial_u h_2 \right) - \left( \partial_u h_1 \cdot \partial_v h_2 \right) \right] \partial_z f \cdot g_2 =0 ,\\
        (\partial_v h_1)(\partial_u h_2)(\partial_x f) - (\partial_u h_1)(\partial_v h_2)(\partial_x f) - (\partial_x h_1)(\partial_u h_2)(\partial_v f) =0.
    \end{aligned}\label{psp:condition}
\end{equation}
Then $P$ is said to be the pseudo-singular point. 
\end{definition}

\begin{figure}[ht]
    \centering
    \begin{subfigure}[t]{0.49\textwidth}
           \includegraphics[width=\linewidth, height=5.5cm]{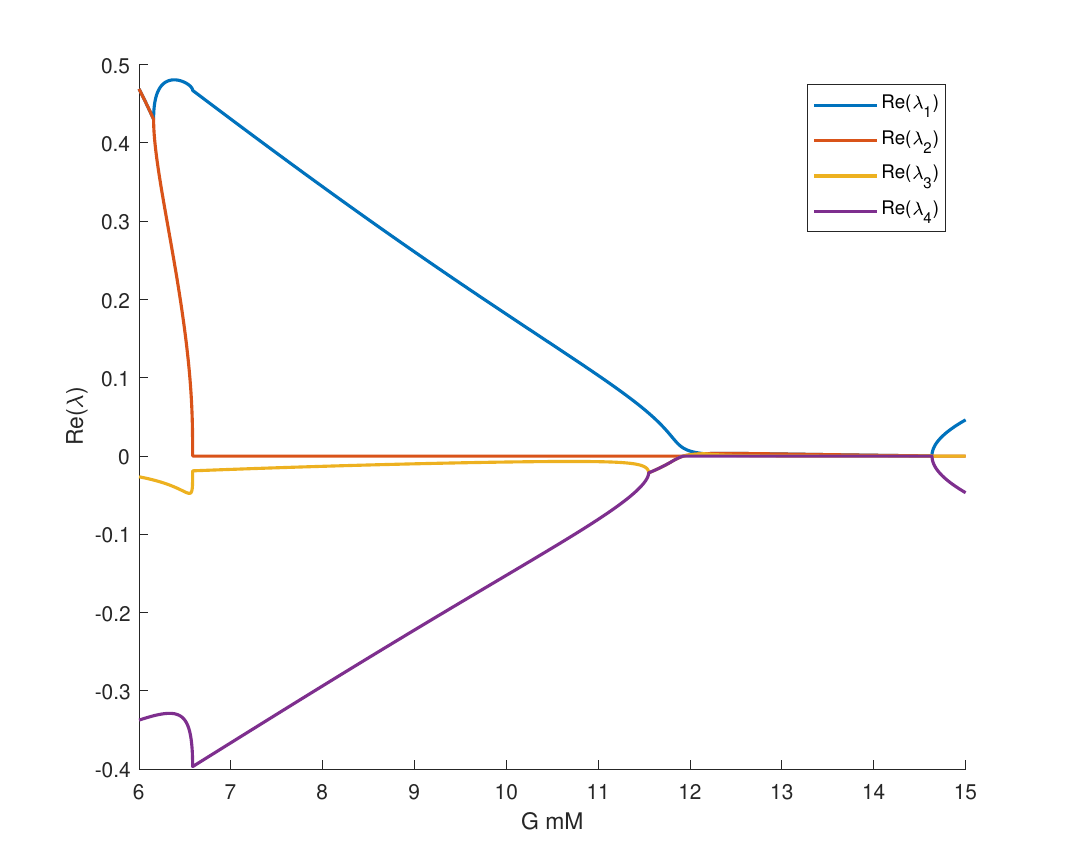}
           \label{fig:Ng1a}
    \end{subfigure}
    \hfill 
    \begin{subfigure}[t]{0.49\textwidth}
        \centering
        \begin{subfigure}[t]{\textwidth}
           \includegraphics[width=\linewidth, height=5.5cm]{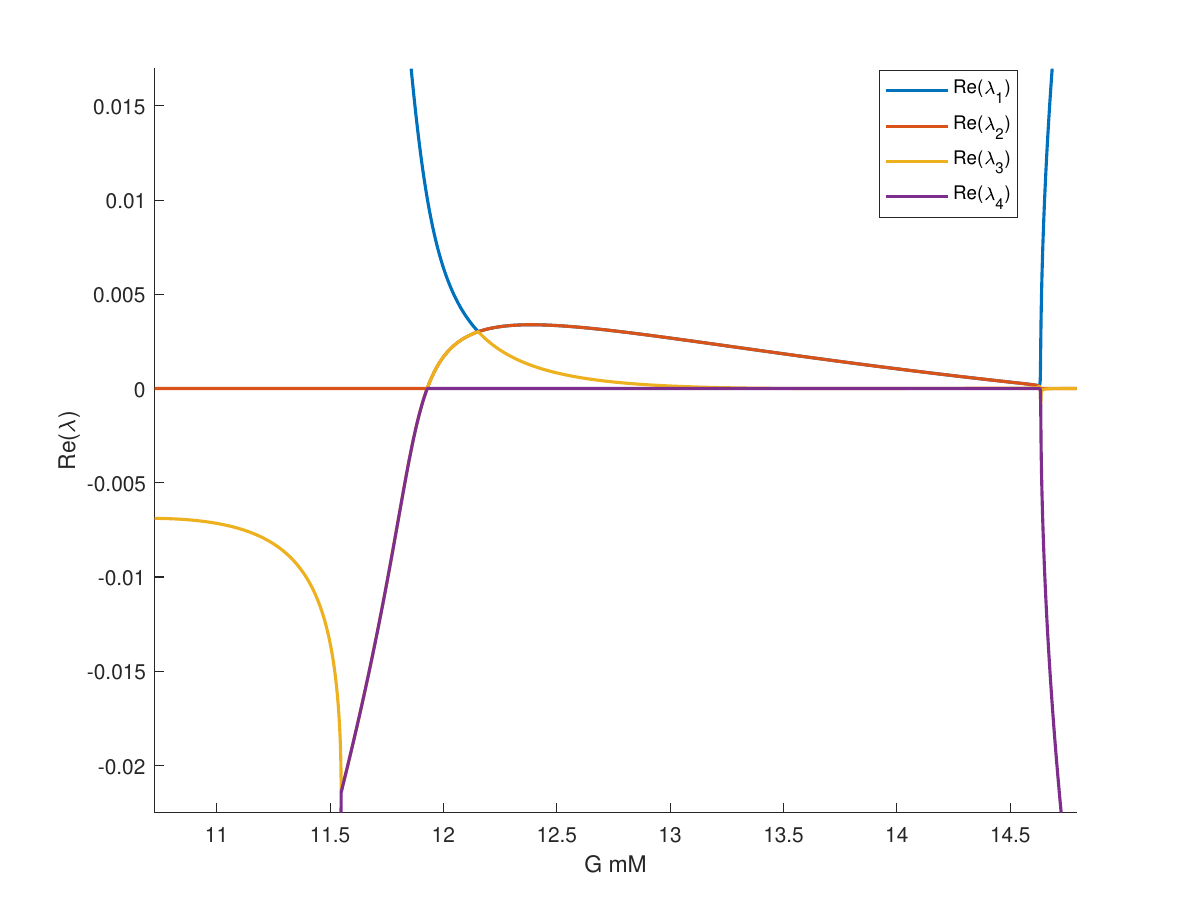}
           \label{fig:Ng1b}
        \end{subfigure}
    \end{subfigure}
    \caption{Eigenvalue at the pseudo singular point of the system (\ref{eq:PSPatzero4D}) for glucose concentration, $G$ ($mM$), in $[6,15]$.}
    \label{fig:Eigen}
\end{figure}

\begin{definition}[Pseudo-singular saddle/node]\label{pseudoTypedefinition}
    Let $p$ be a pseudo-singular point. If the linearized system (\ref{normalized5Dsystem}) at $P$ have two non-zero eigenvalues with opposite sign (plus/minus) then $P$ is called \textit{pseudo singular saddle}, and if all the eigenvalues have same sign then $P$ is called \textit{pseudo singular node} (Figure~\ref{fig:Eigen}).
\end{definition}

\begin{theorem}[Benoit Theorem \cite{benoit1990canards, tchizawa2007generic}]
    If the system (\ref{normalized5Dsystem}) has a pseudo-singular saddle point, then the system has canard solutions. If the system (\ref{normalized5Dsystem}) has a pseudo-singular node point with no resonance, then the system has canard solution (Definition \ref{def:canard}).
\end{theorem}

\begin{assumption}\label{ASSUMPcanardPSP}
    Let $P$ be a pseudo-singular point. At point $P \in \mathcal{C}_1$,
    \begin{equation*}
        \begin{aligned}
            |D_{(v,u)}(h_1, h_2)| =0,\\
            \frac{\partial^2 h_1}{\partial v^2}\neq0,  \ g_2 \neq 0,  \frac{\partial h_1 }{\partial x} \neq 0, \frac{\partial f }{\partial v} \neq 0.
        \end{aligned}
    \end{equation*}
\end{assumption}
\begin{remark}
    Assumption \ref{ASSUMPcanardPSP} is the extension of Benoit’s generic hypothesis \cite{benoit1990canards, ginoux2016canards}, where the non-equality represents the non-degeneracy of the folded singularity.
\end{remark}

Without loss of generality, the following system describes the dynamics in the shifted coordinates, 
\begin{equation}
\begin{aligned}
\frac{dv}{dt} &= -\left( a_1 \hat{X}_\infty(v) (v + v_p - a_2) + \left( a_3 (u + u_p) + a_4 \hat{U}_\infty(x) + a_5 (y + y_p) \right)(v + v_p - a_6) \right):=\hat{h}_1, \\
\frac{du}{dt} &= a_7 \left(\hat{Y}_\infty (v) -u-u_p \right):=\hat{h}_2,\\
\frac{dx}{dt} &= \varepsilon \left( -d_1 \hat{X}_\infty(v) (v + v_p - a_2) - d_2 (x + x_p) + d_3 (z + z_p) \right):=\varepsilon \hat{f}, \\
\frac{dy}{dt} &= \varepsilon \delta \left( \hat{Z}_\infty(x) - (y + y_p) \right):=\varepsilon \delta \hat{g}_1, \\
\frac{dz}{dt} &= \varepsilon \delta \left( k_1 (x + x_p) - k_2 (z + z_p) \right):=\varepsilon \delta \hat{g}_2,
\end{aligned}\label{eq:PSPatzero}
\end{equation}
where \( v, x, y, z \) represent the shifted variables as $v = v_{\text{old}} - v_p$, $u = u_{\text{old}} - u_p$, $x = x_{\text{old}} - x_p$, $y = y_{\text{old}} - y_p$, $z = z_{\text{old}} - z_p$. By the subscript \textit{old}, we mean the original variables of the system (\ref{ScaledMarinelli}. The auxiliary functions are become
\begin{equation*}
\begin{aligned}
\hat{X}_\infty(v) &= \left[ 1 + \exp\left( \frac{v_1 - (v + v_p)}{s_1} \right) \right]^{-1}, \\
\hat{Y}_\infty(v) &= \left[ 1 + \exp\left( \frac{v_2 - (v + v_p)}{s_2} \right) \right]^{-1}, \\
\hat{Z}_\infty(x) &= \left[ 1 + \exp\left( \frac{\frac{G - p_r}{k_r} - (x + x_p)}{s_3} \right) \right]^{-1}, \\
\hat{U}_\infty(x) &= \left[ 1 + \left( \frac{K_d}{x + x_p} \right)^5 \right]^{-1}.
\end{aligned}
\end{equation*}
Now the system (\ref{eq:PSPatzero}) has the pseudo-singular point at the origin.\par
Since $\partial_u\hat{h}_2\neq0$ in the neighborhood of the pseudo singular point, we can express $u=u_p+\hat{Y}_\infty (v)$ on the critical manifold
\begin{equation*}
    \begin{aligned}
        \hat{\mathcal{C}} = \left\{(v,x,y,z) \in \mathbb{R}^4 \ \lvert \ y = -\frac{a_1 \hat{X}_\infty (v) (v + v_p - a_2)}{a_5 (v + v_p - a_6)} - \frac{a_3 \hat{Y}_\infty (v)}{a_5}  - \frac{a_4 \hat{U}_\infty (x)}{a_5}- \frac{a_3 u_p}{a_5} - y_p \right\},
    \end{aligned}
\end{equation*}
and the system (\ref{eq:PSPatzero}) can be written as a $4-$dimensional three-scale systems with $1-$fast variable $v$, $1-$intermediate variable $x$, and $2-$slow variables $y$ and $z$, and reduced system becomes
\begin{equation}
\begin{aligned}
\frac{dv}{dt} &= -\left( a_1 \hat{X}_\infty(v) (v + v_p - a_2) + \left( a_3 \hat{Y}_\infty(v) + a_4 \hat{U}_\infty(x) + a_5 (y + y_p) \right)(v + v_p - a_6) \right):=\hat{h},\\
\frac{dx}{dt} &= \varepsilon \left( -d_1 \hat{X}_\infty(v) (v + v_p - a_2) - d_2 (x + x_p) + d_3 (z + z_p) \right):=\varepsilon \hat{f}, \\
\frac{dy}{dt} &= \varepsilon \delta \left( \hat{Z}_\infty(x) - (y + y_p) \right):=\varepsilon \delta \hat{g}_1, \\
\frac{dz}{dt} &= \varepsilon \delta \left( k_1 (x + x_p) - k_2 (z + z_p) \right):=\varepsilon \delta \hat{g}_2.
\end{aligned}\label{eq:PSPatzero4D}
\end{equation}

\subsection{Normal Form Reduction}
Consider that Assumption (\ref{ASSUMPcanardPSP}), and condition (\ref{psp:condition}) satisfied at $(0,0,0,0)$. Then, by rescaling $v,x,y,z$ as $V= \frac{\hat{h}_{vv}}{2}v$, $X= -\frac{\hat{h}_{vv} \hat{h}_x}{2}x$, $Y= \frac{\hat{h}_{vv} \hat{h}_y}{2}y$ and $Z=\frac{1}{\hat{g}_2(0)}z$, the system (\ref{eq:PSPatzero}) can be written in the form with new variables $V,X,Y$ and $Z$,
\begin{equation}\label{NORMALfORMbETA}
\begin{aligned}
    \frac{dV}{dt} &= -X H_1(V, X, Y) + Y H_2(V, X, Y) + V^2 H_3(V, X, Y),\\
\frac{dX}{dt} &= \varepsilon \left( b_1 V H_5(V, X, Z) + b_2 Z H_6(V, X, Z) + b_3 X H_7(V, X, Z) \right),\\
  \frac{dY}{dt} &= \varepsilon \delta \left( -\lambda H_8(V, X, Y, Z) - b_4 X H_9(V, X, Y, Z) -b_5 Y \right),\\
\frac{dZ}{dt} &= \varepsilon \delta H_{10}(V, X, Y, Z),
\end{aligned}
\end{equation}
where, we only presented up to third-order monomials, 
\allowdisplaybreaks
    \begin{align*}
    H_1(V, X, Y) &= 1 - H_{XX} X - H_{VX} V - H_{XXX} X^2 +\mathcal{O}(X^3),\\
        H_2(V, X, Y) &= 1  + H_{VY} V,\\
        H_3(V, X, Y) &= 1 + H_{VVV} V  + \mathcal{O}(V^2),\\
        H_5(V, X, Z) &= 1 + \frac{F_{VV}}{F_V} V + \frac{F_{VVV}}{F_V} V^2 + \mathcal{O}(V^3),\\
        H_6(V, X, Z) &= 1,\\
        H_7(V, X, Z) &= 1,\\
        H_8(V, X, Y, Z) &= 1,\\
        H_9(V, X, Y, Z) &= 1 + \frac{G_{1XX}}{G_{1X}} X + \frac{G_{1XXX}}{G_{1X}} X^2 + \mathcal{O}(X^3),\\
        H_{10}(V, X, Y, Z) &= 1 + G_{2X} X + G_{2Z} Z,\\
        b_1 &= F_V, \\
        b_2 &= F_Z, \\
        b_3 &= F_X, \\
        b_4 &= -G_{1X}, \\
        b_5 &= G_{1Y}, \\
        \lambda &= -G_{10}.
    \end{align*}
\allowdisplaybreaks[0]

With the relation between the reduced $4D$ system (\ref{eq:PSPatzero}) and reduced normal form (\ref{NORMALfORMbETA}), the coefficients can be written as
\begin{equation*}
    \begin{aligned}
    \renewcommand{\arraystretch}{1.5}
\begin{tabular}{ll@{\hspace{2cm}}ll}
            $H_{XX}$ &= $\dfrac{\hat{h}_{xx}}{\hat{h}_{vv} \hat{h}_x^2}$, & $G_{10}$ &= $\dfrac{\hat{h}_{vv} \hat{h}_y \hat{g}_{1}(0)}{2}$, \\
            $H_{VX}$ &= $-\dfrac{2 \hat{h}_{vx}}{\hat{h}_{vv} \hat{h}_x}$, & $G_{1X}$ &= $-\dfrac{\hat{h}_y \hat{g}_{1x}}{\hat{h}_x}$, \\
            $H_{VY}$ &= $\dfrac{2 \hat{h}_{vy}}{\hat{h}_{vv} \hat{h}_y}$, & $G_{1Y}$ &= $\hat{g}_{1y}$, \\
            $H_{VVV}$ &= $\dfrac{2 \hat{h}_{vvv}}{3 \hat{h}_{vv}^2}$, & $G_{1XX}$ &= $\dfrac{\hat{h}_y \hat{g}_{1xx}}{\hat{h}_{vv} \hat{h}_x^2}$, \\
            $H_{XXX}$ &= $-\dfrac{2 \hat{h}_{xxx}}{3 \hat{h}_{vv}^2 \hat{h}_x^3}$, & $G_{1XXX}$ &= $-\dfrac{2 \hat{h}_y \hat{g}_{1xxx}}{3 \hat{h}_{vv}^2 \hat{h}_x^3}$, \\
            $F_V$ &= $-\hat{h}_x \hat{f}_v$, & $G_{2X}$ &= $-\dfrac{2 \hat{g}_{2x}}{\hat{h}_{vv} \hat{h}_x \hat{g}_{2}(0)}$, \\
            $F_X$ &= $\hat{f}_x$, & $G_{2Z}$ & $= \hat{g}_{2z}$, \\
            $F_Z$ &= $-\dfrac{2\hat{f}_z \hat{g}_{2}(0)}{\hat{h}_{vv} \hat{h}_x}$, & $F_{VV}$ &= $-\dfrac{\hat{h}_x \hat{f}_{vv}}{\hat{h}_{vv}}$, \\
            $F_{VVV}$ &= $-\dfrac{2 \hat{h}_x f_{vvv}}{3 \hat{h}_{vv}^2}$. & & \\
        \end{tabular}
    \end{aligned}
\end{equation*}
Here, the subscript with $v,x,y,z$ means the partial derivative, and $\mathcal{O}_3$ means 3rd and higher order monomial. All partial derivatives are evaluated at the (shifted) pseudo-singular point $0$. The coefficient values are presented in Table~\ref{coeffficientBETAnormal}.

\begin{table}[ht]
\centering
\begin{tabular}{l l l l}
\toprule
\textbf{Coefficient} & \textbf{Value} & \textbf{Coefficient} & \textbf{Value} \\
\midrule
$H_{XX}$     & $-3.42120 \times 10^{6}$  & $F_{VV}$     & $2.93710 \times 10^{-3}$ \\
$H_{VX}$     & $-4.99657 \times 10^{1}$  & $F_{VVV}$    & $3.76990 \times 10^{-2}$ \\
$H_{VY}$     & $4.99657 \times 10^{1}$   & $G_{10}$     & $-1.00000 \times 10^{-4}$ \\
$H_{VVV}$    & $1.15445 \times 10^{1}$   & $G_{1X}$     & $-7.66910 \times 10^{2}$ \\
$H_{XXX}$    & $3.58420 \times 10^{13}$  & $G_{1Y}$     & $-1.00000 \times 10^{0}$ \\
$F_V$        & $1.26600 \times 10^{-4}$  & $G_{1XX}$    & $2.20800 \times 10^{7}$ \\
$F_X$        & $-1.60000 \times 10^{1}$  & $G_{1XXX}$   & $1.68630 \times 10^{14}$ \\
$F_Z$        & $-5.37860 \times 10^{-6}$ & $G_{2X}$     & $-1.44480 \times 10^{7}$ \\
              &        & $G_{2Z}$     & $-7.46300 \times 10^{0}$ \\
\bottomrule
\end{tabular}
\caption{Coefficient values for (\ref{NORMALfORMbETA}) at $(v_p, x_p, y_p, z_p)= (-60.4, 0.094, 0.467, 84.539)$  which is pseudo-singular saddle (Definition~\ref{pseudoTypedefinition}).}
\label{coeffficientBETAnormal}
\end{table}

\subsection{Blow-Up at the Pseudo-Singular Point}
For the singular perturbation problems, Geometric Singular Perturbation Theory (GSPT) \cite{fenichel1979geometric, kuehn2015multiple, guckenheimer2009computing}is a powerful tool to analyze the dynamics of the system by dissecting the fast and slow dynamics on the critical manifold when it is hyperbolic (real parts of all the eigenvalues along the fast direction are non-zero). In system (\ref{NORMALfORMbETA}), the origin is on the fold curve of the critical manifold $\hat{\mathcal{C}}$, which is a non-hyperbolic equilibrium (has zero eigenvalue in the fast direction) of the system, where we cannot apply this theory. To analyze the dynamics in the neighborhood of the non-hyperbolic point, first, we have to desingularize the system at that point. The blow-up technique is the right tool to desingularize the non-hyperbolic points \cite{desroches2012mixed, szmolyan2001canards,chiba2011periodic, jardon2019survey}. 
The blow-up methods consist of the local change of variables and are applied to desingularize the flows near the fold point. 
\begin{definition}
    The blow-up transformation $\Phi$, defined by 
\begin{equation*}
\begin{aligned}
        \Phi\ : \ \mathbb{S}^7 \times \mathbb{R}^+ &\to \mathbb{R}^7,\\
    (\bar{v}, \bar{x}, \bar{y}, \bar{z}, \bar{\varepsilon}, \bar{\delta}, \bar{\lambda}, \bar{r})  &\to (\bar{r}\bar{v}, \bar{r}^2\bar{x}, \bar{r}^4\bar{y}, \bar{r}^2\bar{z}, \bar{r}^2\bar{\varepsilon}, \bar{r}\bar{\delta}, \bar{r}^2\bar{\lambda})=(v,x,y,z,\varepsilon,\delta,\lambda).
\end{aligned}
\end{equation*}
where $\mathbb{S}^7$ is a $7-$sphere.
\end{definition}

\begin{definition}[Directional Blow-Up]
    For system (\ref{NORMALfORMbETA}), directional blow-ups $\Phi_i,\ i=1,\cdots,14$ are obtained by settings one blown-up variables on $\mathbb{S}^7$ equal to $\pm 1$ in the definition of $\Phi$. The directional charts $K_i$ are defined such that $\Phi(\mathbb{S}^7 \times \mathbb{R}^+)=\Phi_i(K_i(\mathbb{S}^7 \times \mathbb{R}^+))$.
\end{definition}

In the singular perturbation problems, the chart along the direction where $\bar{\varepsilon}=1$ is the most important chart, called the central chart, or classical chart. In this paper, we will only present here the numerical results of the central chart, say, $K_3$. 
\begin{align*}
    V=r_3 v_3,\quad X= r_3^2 x_3, \quad Y=r_3^4 y_3, \quad Z = r_3^2 z_3, \quad \varepsilon = r_3^2 , \quad \delta = r_3 \delta_3 , \quad \lambda= r_3^2 \lambda_3,
\end{align*}
where $\varepsilon_3=1$.

Using this coordinate, in $K_3$ the system (\ref{NORMALfORMbETA}) with additional differential equations
\begin{align*}
    \frac{d\varepsilon}{dt} =0,\ \quad \frac{d\delta}{dt} =0, \quad \frac{d\lambda}{dt} =0,
\end{align*}
becomes
\begin{equation}
    \begin{aligned}
        \frac{dv_3}{dt_3}&=-x_3 + v_3^2 + r_3 ( r_3 y_3 + r_3 H_{XX} x_3^2 + H_{VX} v_3 x_3 + r_3^2 H_{VY} v_3y_3 + H_{VVV} v_3^3 \\  & \quad + r_3^3 H_{XXX} x_3^3 + \mathcal{O}(r_3, v_3^4)+ \mathcal{O}(r_3^5,x_3^4))),\\
        \frac{dx_3}{dt_3}&= b_1 v_3 + r_3 (b_2 z_3 + b_3 x_3 + F_{VV} v_3^2 + r_3 F_{VVV} v_3^3 + \mathcal{O}(r_3^2, v_3^4)),\\
        \frac{dy_3}{dt_3} &= \delta_3 \left[-\lambda_3-b_4 x_3 + r_3^2 ( -b_5 y_3 + G_{1XX} x_3^2 + r_3^2 G_{1XXX} x_3^3 + \mathcal{O}(r_3^4, x_3^4))\right],\\
        \frac{dz_3}{dt_3} &= \delta_3 \left[1+r_3^2 (G_{2X} x_3 + G_{2Z} z_3) \right],
    \end{aligned}\label{K3systembetacell}
\end{equation}
and $\frac{dr_3}{dt_3}=0, \ \frac{d\delta_3}{dt_3}=0, \ \frac{d\lambda_3}{dt_3}=0$. Here, $t_3=r_3t = \sqrt{\varepsilon}t$.

\begin{remark}
This transformation is a $\varepsilon$ dependent rescaling of the variables and the parameters, we are interested in their influence on the dynamics of the system (\ref{NORMALfORMbETA}). The exponent of $r$ or $r_3$ is called the \textit{weight}. A theory for obtaining the weight is associated with the Newton Diagram \cite{chiba2013first, chiba2024weights}.
\end{remark}

\begin{remark}
    The relation between the original system (\ref{eq:PSPatzero}) and the $K_3$ coordinate, can be numerically (using the parameters in Table~\ref{tab:ScalePara} and Table \ref{coeffficientBETAnormal}) defined as
    \begin{align*}
        v = 13.838 v_3, \quad x = 2.3061 \times 10^{3} x_3, \quad y = -7.414 \times 10^{-5} y_3, \quad z = -0.253 z_3, \quad t = 51.6191 t_3, 
    \end{align*}
since $r_3=\sqrt{\varepsilon}$.
\end{remark}

In the invariant manifold $r_3=0$, the system (\ref{K3systembetacell}) has special solution
\[\gamma(t_3):=(v_3(t_3),x_3(t_3), y_3(t_3),z_3(t_3))=(\frac{b_1}{2}t_3, \frac{b_1^2}{4}t_3^2-\frac{b_1}{2}, \delta_3((\frac{b_1 b_4}{2}-\lambda_3)t_3-\frac{b_1^2 b_4}{12} t^3_3),\delta_3 t_3+z_{30}),\]
where $z_{30}$ is the initial value of $z_3$.

\begin{figure}[ht]
    \centering
    \hspace*{-1.1cm} 
    \includegraphics[width=0.9\textwidth]{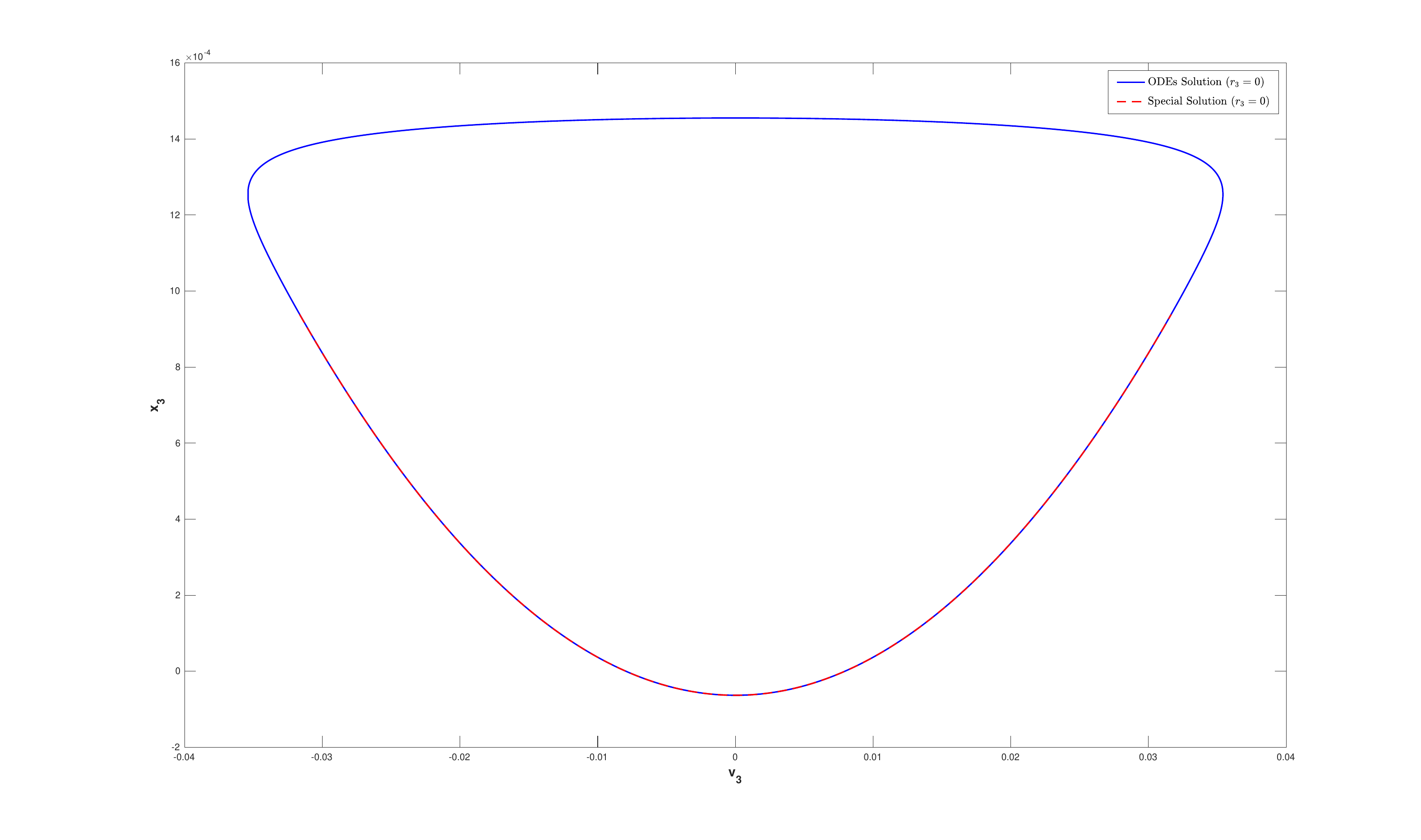}
    \caption{In the invariant manifold $r_3=0$, the special solution (red) and the solution of the system (\ref{K3systembetacell}) (blue) in $(v_3, x_3)$.}
    \label{special_solution_G_6815_r_3_0}
\end{figure}

The unperturbed solution of the system (\ref{K3systembetacell}) is the separatrix along the special solution, which separates closed orbits from unbounded orbits (Figure~\ref{special_solution_G_6815_r_3_0}).

\begin{figure}[ht]
    \centering
    \begin{subfigure}[b]{0.48\textwidth}
        \centering
        \includegraphics[width=\textwidth]{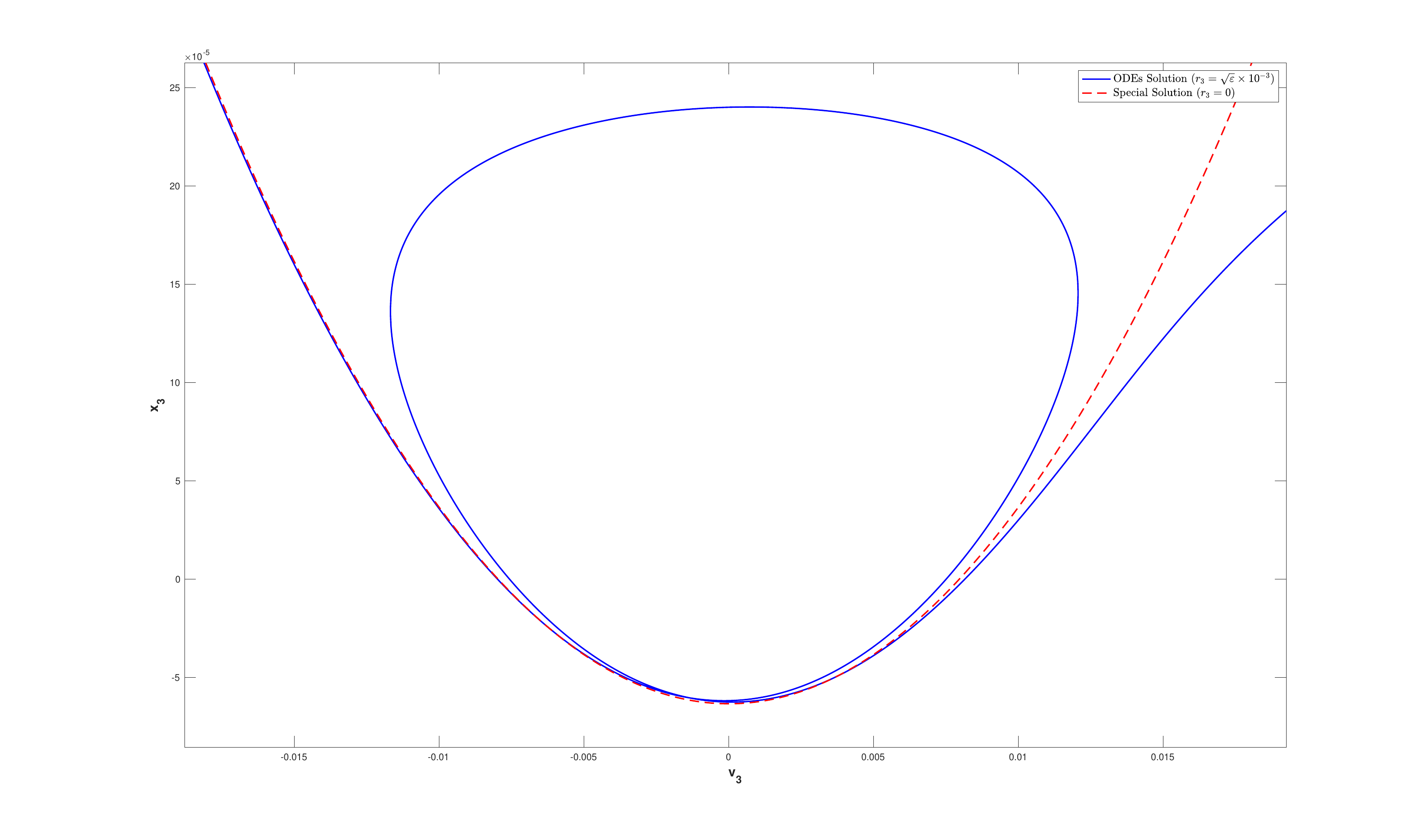}
        \caption{$r_3=\sqrt{\varepsilon}\times 10^{-3}$}
        \label{special_solution_G_6815_r_3_sqrt_eps_time10e-3}
    \end{subfigure}
    \hfill
    \begin{subfigure}[b]{0.48\textwidth}
        \centering
        \includegraphics[width=\textwidth]{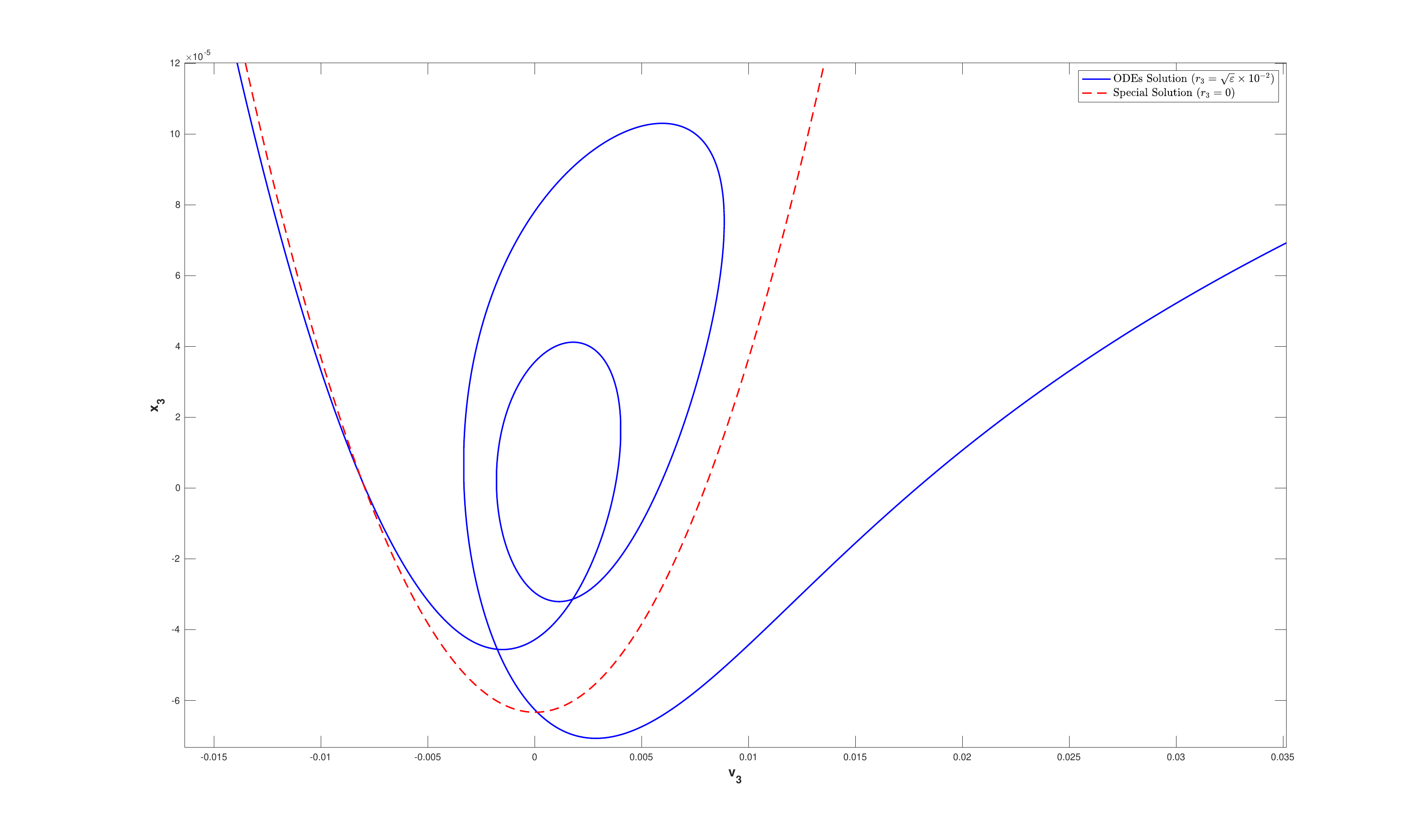}
        \caption{$r_3=\sqrt{\varepsilon}\times 10^{-2}$}
        \label{special_solution_G_6815_r_3_sqrt_eps_time10e-2}
    \end{subfigure}

    \vspace{0.5cm}

    \begin{subfigure}[b]{0.48\textwidth}
        \centering
        \includegraphics[width=\textwidth]{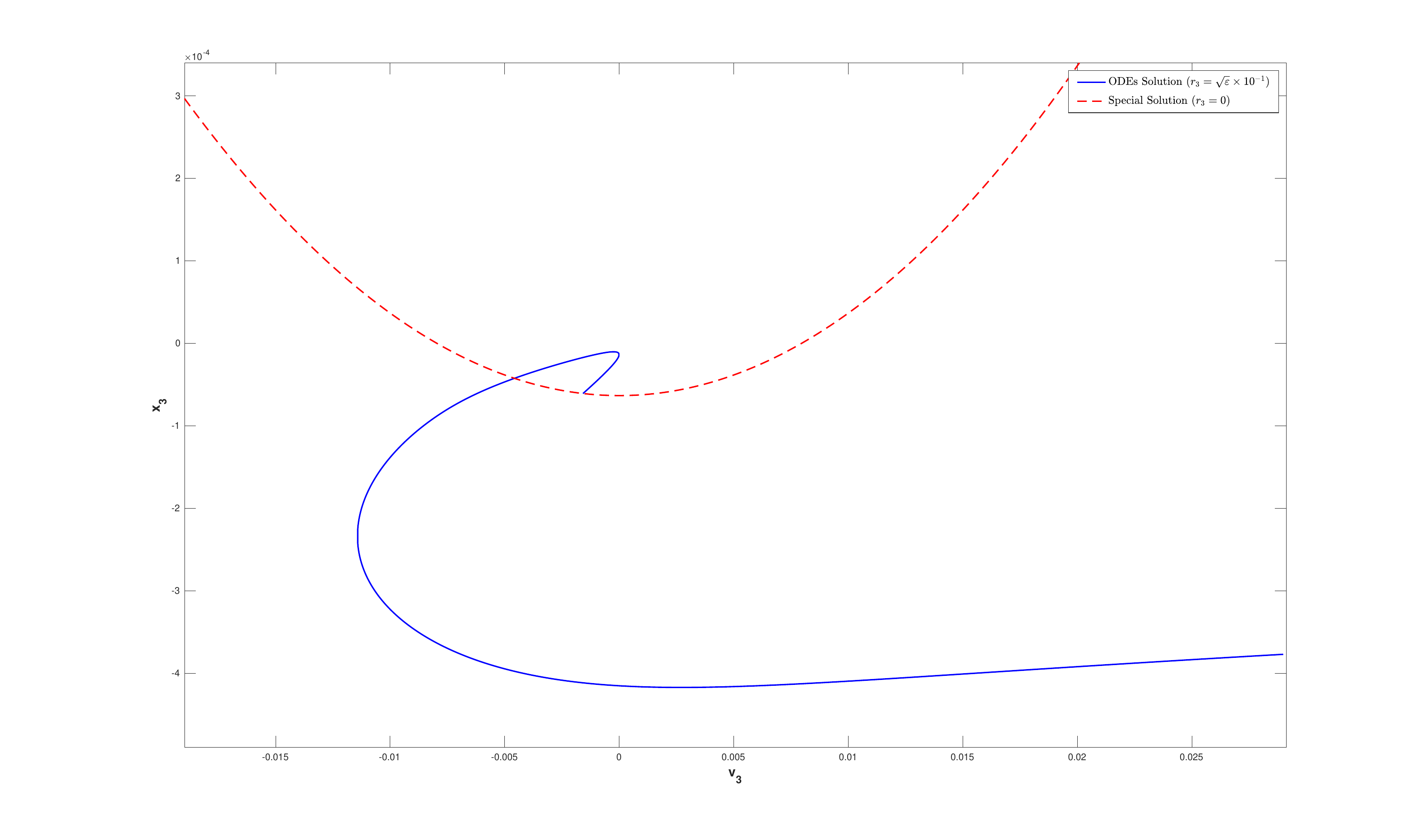}
        \caption{$r_3=\sqrt{\varepsilon}\times 10^{-1}$}
        \label{special_solution_G_6815_r_3_sqrt_eps_time10e-1}
    \end{subfigure}
    \hfill
    \begin{subfigure}[b]{0.48\textwidth}
        \centering
        \includegraphics[width=\textwidth]{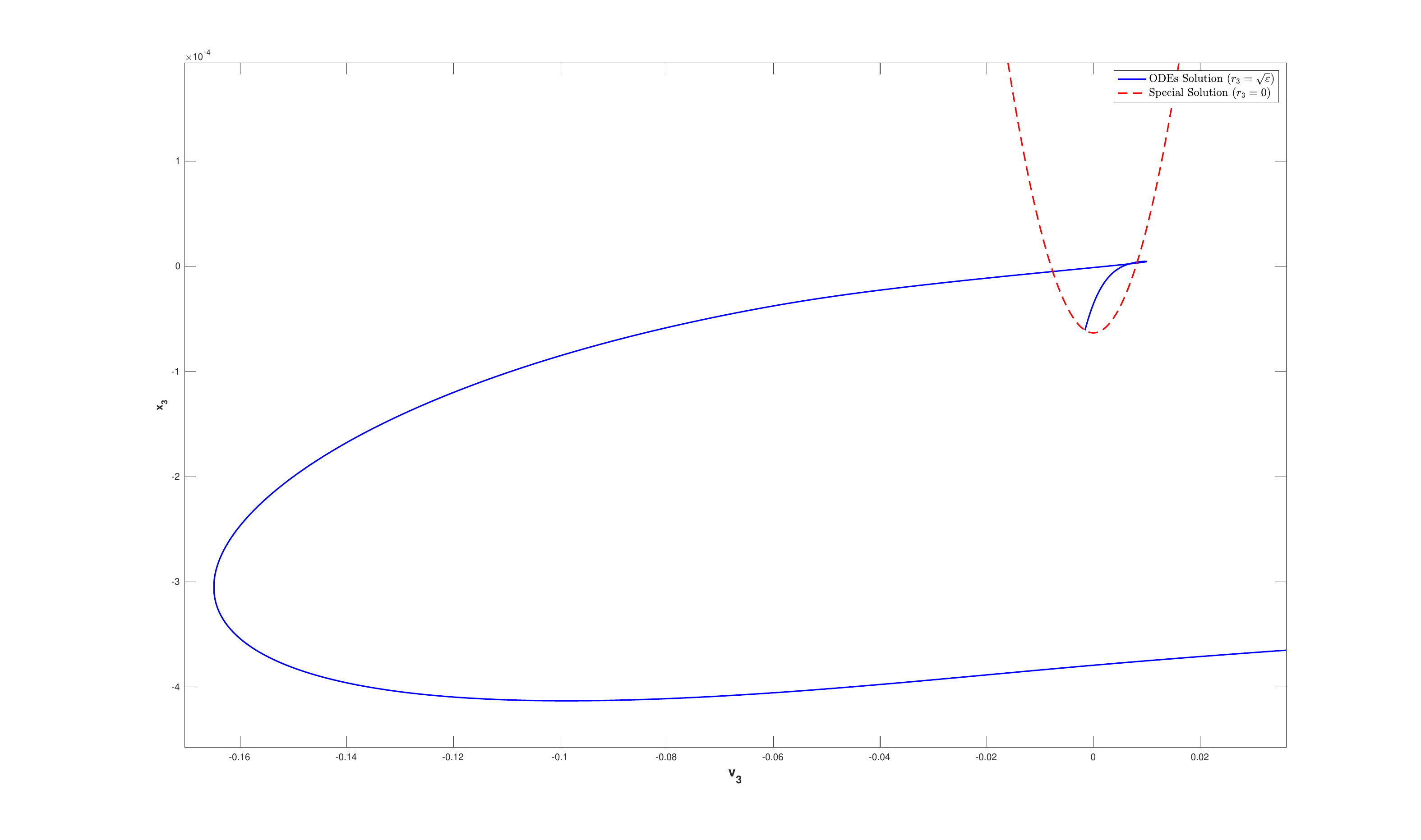}
        \caption{$r_3=\sqrt{\varepsilon}$}
        \label{special_solution_G_6815_r_3_sqrt_eps}
    \end{subfigure}

    \caption{Solution of perturbed system (\ref{K3systembetacell}) (blue) compared to the special solution (red) in $(v_3,x_3)$.}
    \label{fig:k3r3}
\end{figure}

As the $r_3$ grows from $0$ to $\sqrt{\varepsilon}$, the perturbed solution creates small oscillations near the focus of the separatrix before growing unboundedly (Figure~\ref{special_solution_G_6815_r_3_sqrt_eps_time10e-3}, \ref{special_solution_G_6815_r_3_sqrt_eps_time10e-2}). For sufficiently large $r_3$, this solution gets negative feedback because of the solwest variable dynamics, before growing in the fast ($v_3$) direction unboundedly (Figure~\ref{special_solution_G_6815_r_3_sqrt_eps_time10e-1}, \ref{special_solution_G_6815_r_3_sqrt_eps}).

\section{Numerical Results}\label{sec:numerical}
\subsection{Limit Cycle and \Poincare Return Map}
\begin{figure}[ht]
    \centering
    \hspace*{-1.1cm}
    \includegraphics[width=0.6\textwidth]{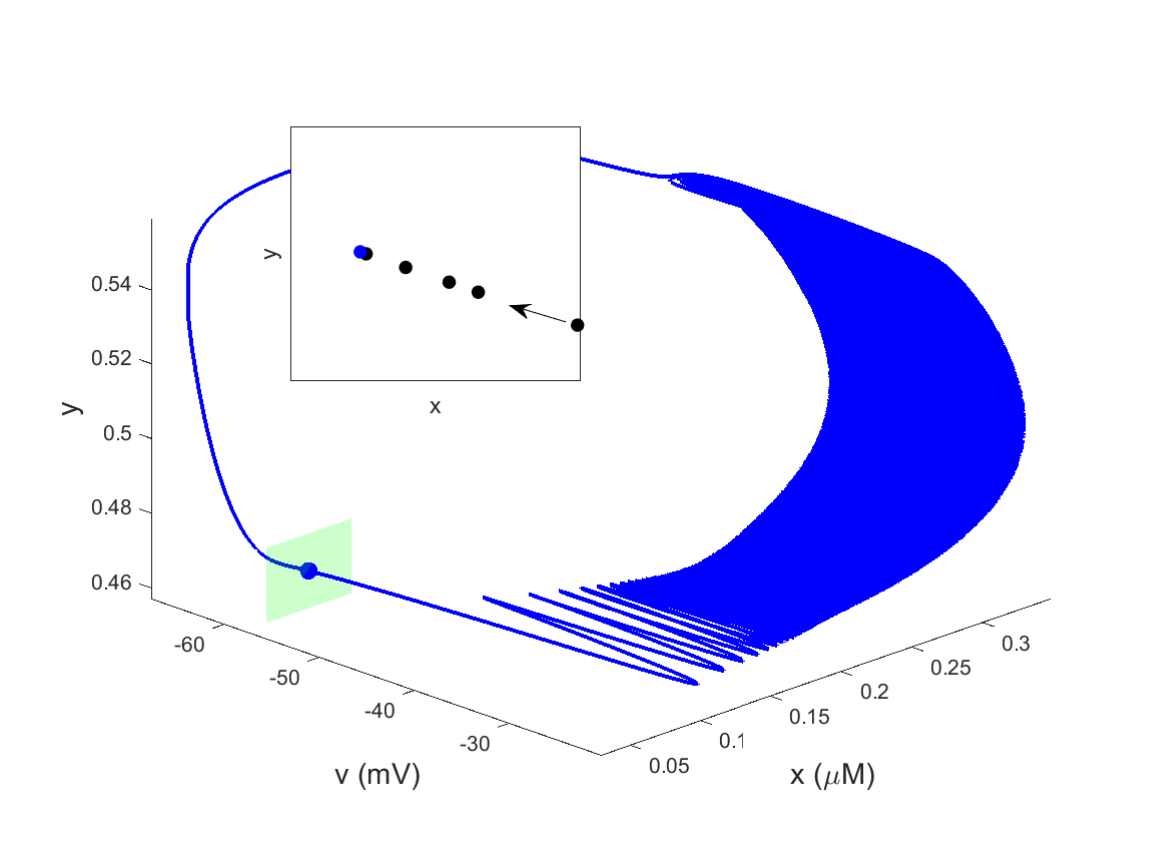}
    \caption{Limit cycle in $(v, x, y)-$space (blue curve) with \Poincare Section (green surface) at $v=-59\ mV$ for $G = 8\ mM$, showing a stable fixed point on the section. Inset: \Poincare map in $(x, y)-$plane, with arrows indicating convergence to the fixed point.}
    \label{FIG:limitCYCLE_G_8}
\end{figure}

Figure~\ref{FIG:limitCYCLE_G_8} illustrates the limit cycle of the model (\ref{ScaledMarinelli}) in the ($ v, x, y$)-space for $G = 8\ mM$, where $v $ is the membrane potential (mV), $x$ is cytosolic calcium ($\mu M$), and $y$ is a slow variable (ADP/ATP ratio) modulating bursting. The blue curve represents the periodic orbit, capturing the slow oscillatory dynamics with distinct bursting phases. The green surface at $v = -59$ is the Poincaré section, intersecting the limit cycle at points corresponding to burst initiation. The inset displays the Poincaré map in the $(x, y)$-plane, where each point represents the state at a section crossing. The map shows convergence to a stable fixed point (e.g., at $x \approx 0.073$,$ y \approx 0.467$), with arrows indicating the iterative approach to this point. This fixed point confirms the stability of the limit cycle, reflecting a consistent burst initiation pattern.\par

\begin{figure}[ht]
    \centering
    \includegraphics[width=0.6\textwidth]{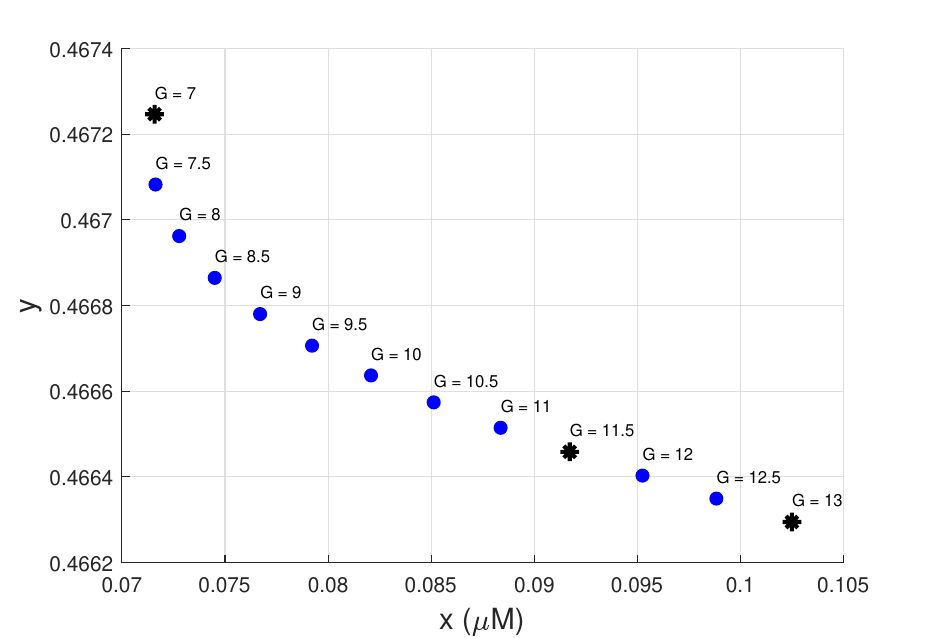}
    \caption{Fixed points of the \Poincare map in the ($x, y$)-plane for $G \in [7, 13]$ mM. Blue circles indicate stable limit cycles, and black stars mark where the trajectory remains close to the limit cycles over the time windows but does not monotonically converge to the fixed point on \Poincare section. Labels denote $G$ values.}
    \label{fig:poincare_map}
\end{figure}

\begin{table}[ht]
    \centering
    \begin{tabular}{cc}
        \begin{minipage}{0.45\textwidth}
            \centering
            \includegraphics[width=1.1\linewidth]{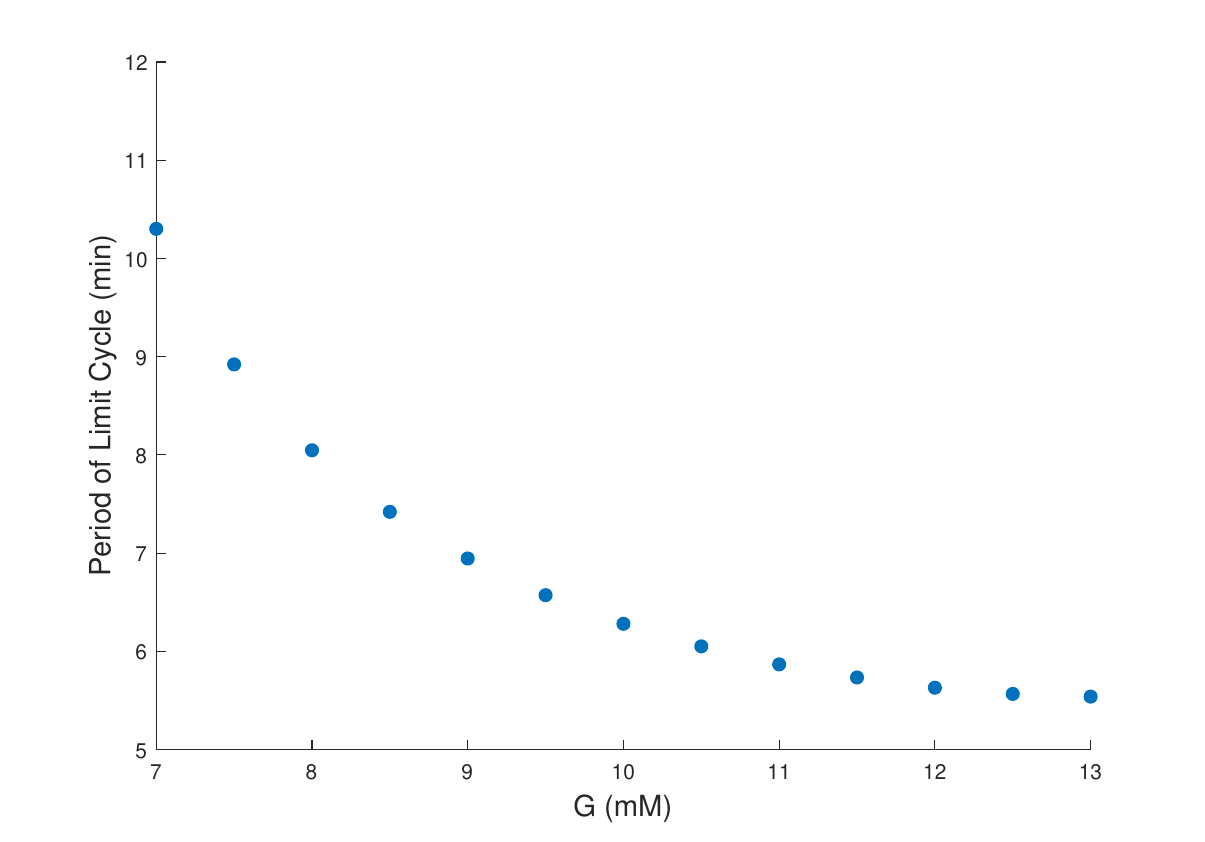}
            \captionof{figure}{Period of the limit cycle (minutes) vs. $G \in [7,13]$ (mM).}
            \label{fig:period_vs_G}
        \end{minipage}
        &
        \begin{minipage}{0.45\textwidth}
            \centering
            \includegraphics[width=\linewidth]{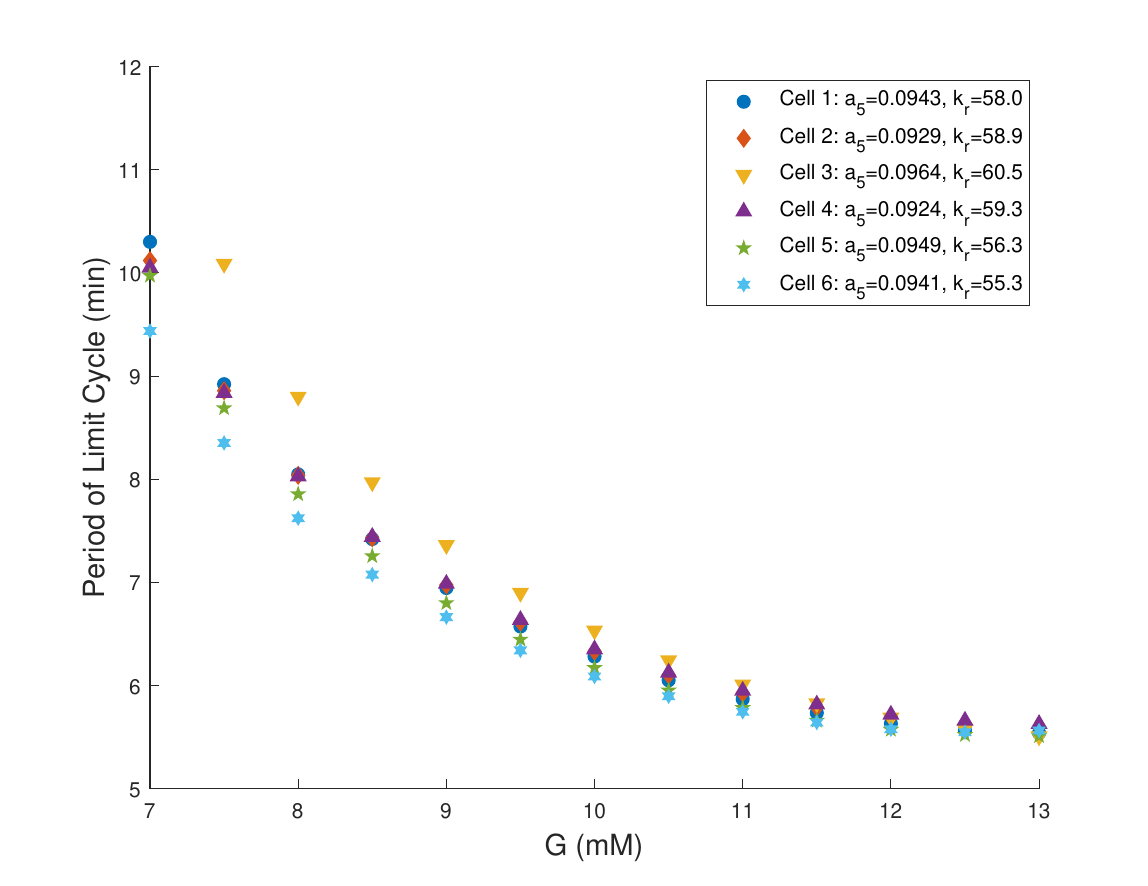}
            \captionof{figure}{Period of the limit cycle (minutes) for 6 cells vs. $G \in [7,13]$ (mM).}
            \label{fig:period_vs_G_6_cell}
        \end{minipage}
    \end{tabular}
\end{table}

The $\beta-$cell model (\ref{ScaledMarinelli}) exhibits slow oscillations with bursting behavior (Figure~\ref{fig:single_beta_burst}) driven by glucose ($G$). A \Poincare section at $v = -59$ mV, where $v$ is the membrane potential, captures burst initiation via upward crossings. Simulations over $t \in [0, 80]$ ($min$) for $G \in [7, 13]$ mM reveal $8$ to $15$ crossings per $G$, indicating varying burst frequencies. The \Poincare map records state variables $(x,y)$ at each crossing, where $x$ is cytosolic calcium ($\mu$M) and $y$ is the ADP/ATP ratio. For each $G$, fixed points are computed using the last 6 crossings, with stability assessed via the maximum distance: $max_{dist} = \max_{n}\left(\sqrt{(x_{n+1} - x_{n})^2 + (y_{n+1} - y_{n})^2}\right)$, $n=1,\cdots,5$, using a tolerance of $1 \times 10^{-6}$, and monotonous decreasing condition of this distance (contraction) to distinguish stable limit cycle.

Figure~\ref{fig:poincare_map} shows the fixed points of the \Poincare map in the ($x, y$)-plane for $G$ from 7 to 13 mM. Blue circles represent stable limit cycles ($max_{dist} < 1 \times 10^{-6}$), indicating consistent burst initiation (e.g., $G = 7.5$ to 10.5, and 12 to 12.5, with $max_{dist}$ from $5.30 \times 10^{-8}$ to $2.81 \times 10^{-10}$). For increasing value of $G$, fixed points on the \Poincare section, shift downward in $y$, while $x$ goes around 0.07 to 0.10 $\mu$M, reflecting decreasing quiescent phase.
Figure~\ref{fig:period_vs_G} plots the period of the limit cycle (ms) against $G$ (mM), with a color gradient from blue ($G = 7$) to yellow ($G = 13$). The period decreases from approximately 6.2 ms at $G = 7$ to 3.2 ms at $G = 13$, indicating faster bursting as glucose increases. This trend aligns with the increasing number of crossings (8 to 15), reflecting shorter inter-burst intervals at higher $G$. Stable limit cycles ($G = 7.5$ to 10.5) show a smooth period decrease, while at $G = 7, 11.5, 13$ exhibit slight deviations. This analysis elucidates glucose-dependent dynamics in $\beta-$cell bursting.

\subsubsection*{Impact of Glucose and Heterogeneity}
For uncoupled heterogeneous Pancreatic $\beta-$cells, the periods of the limit cycle at the lower glucose concentration of the range $[7,13]$ have significant differences, in minutes. As the glucose concentration level increases in that range, the period of the limit cycle becomes nearly unnoticeable (Figure~\ref{fig:period_vs_G_6_cell}). Here, by the heterogeneity, we mean the variation in the $a_5$ values, which corresponds to the maximal conductance for $ATP-$sensitive potassium ($K^+$) channel, usually in a biophysical model denoted as $g_{K(ATP)}$. Besides $a_5$, we also choose different values for $k_r$ in system (\ref{ScaledMarinelli}). Since we are only considering a three-time-scale system without the $G$ evolution equation, by choosing $k_r$ at different values, we are artificially reflecting the situation, in the case where the cells have different glucose sensitivity/threshold or different glucose uptake ability.

\subsection{Synchronization in the Local Network of Beta-Cells (Slow Burster)}
The Marinelli et al.\ model \cite{marinelli2022oscillations} exhibits slow oscillations with a period of $343$ seconds ($5.7$ minutes), characterized by $3$ complete bursts in $20$ minutes (Figure \ref{fig:single_beta_burst}), likely reflecting compound bursting. Synchronization of these oscillations across $\beta-$cells within an islet, achieved through electrical coupling, is crucial for coordinated insulin secretion. This coupling is modeled as '$\frac{k}{N}\sum_{j} (v_j - v_i)$', where $v_i$ and $v_j$ are membrane potentials of coupled cells, and $k$ is the coupling coefficient proportional to gap junction conductance \cite{sherman1991model}. More precisely, the coupling strength in system (\ref{nbeta}) is defined by $k=\frac{g_c}{c_v}$, where $g_c$ is the gap junction conductance, and $c_v$ is the membrane capacitance. The range of coupling strength $k$, its dependence on glucose ($G$), and the mathematical study of synchronization are pivotal for understanding the network dynamics, particularly for the initiation of bursting.

Synchronized activity between Pancreatic $\beta-$cells ensures unified insulin pulses. These pulses optimize glucose regulation by enhancing insulin sensitivity. The coupling term drives diffusive current, synchronizing membrane potentials via gap junctions. For slow oscillations, this coordination maximizes calcium influx and insulin release. In bursting, synchronization can involve spike-to-spike alignment, initiation of bursting, or end-of-burst timing. Here, synchronizing the initiation of bursting ensures that all $\beta-$cells start their active phases simultaneously, amplifying the islet’s insulin pulse \cite{sherman1991model}. Desynchronization, as in Type-2 diabetes, disrupts this rhythm, impairing glucose control.

The coupling coefficient $k$, in $ms^{-1}$ ($\left[\frac{pS}{fF}\right] = [ms^{-1}]$), determines synchronization efficacy. Sherman et al.\ \cite{sherman1991model} reported that at $g_c \approx 10-175 \, \text{pS}$ (for our studied model, $k \approx 0.0019-0.033\ (ms^{-1})$, since we rescaled the write hand side of $\frac{dv}{dt}$ equation by capacitance, $c_v = 5300\ fF$), burst initiation is almost synchronized. In our numerical result, $k$ value in the range $[0.005, 0.045]$ $ms^{-1}$ which implies $g_c \in[30,240]$ $pS$ for $G\in[7,13]\ mM$. As $G$ increases from $7 $ to $13\ mM$, coupling rises due to glucose-driven depolarization, increasing conductance.

The glucose concentration ($G$) modulates $\beta$-cell oscillatory dynamics by altering glycolytic flux and KATP channel activity, which may influence the coupling strength required for synchronization. In our study, increasing glucose from $G \in [7, 13] \ mM$ reduces the bursting period from approximately 10.3 minutes to 5.5 minutes (period of limit cycle), indicating faster glycolytic oscillations and KATP conductance dynamics. Limited literature directly addresses how glucose levels affect coupling strength requirements.

\begin{figure}[ht]
    \centering
    \includegraphics[width=0.8\textwidth]{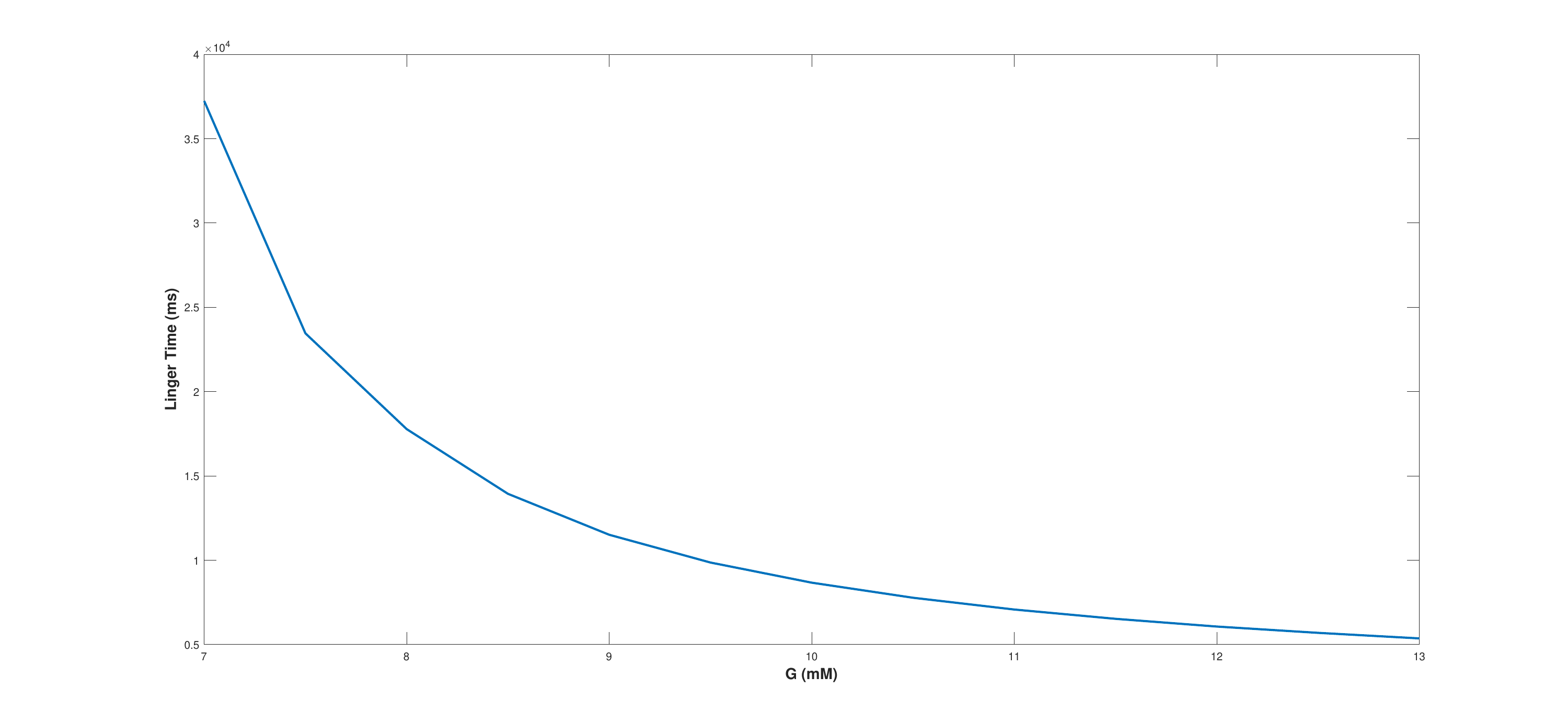} 
    \caption{Linger time ($t_{linger}\ (ms)$) as a function of glucose concentration ($G$).}
    \label{fig:linger_time_G}
\end{figure}

Figure \ref{fig:linger_time_G} shows the minimum linger time as a function of glucose concentration ($G$) for a single pancreatic $\beta-$cell. The linger time ($t_{linger}$) is the duration the trajectory spends in the neighborhood of the pseudo-singular point, corresponding to the quiescent phase of the bursting cycle. The x-axis represents (G) in millimolar ($mM$), ranging from $7$ to $13\ mM$, covering normoglycemic to hyperglycemic conditions. 
The linger time ($t_{linger}$) plot against glucose concentration ($G$) displays a decreasing trend, with the linger time dropping from approximately $35,000$ ms at $G = 7\ mM$ to $5,000\ ms$ at $G = 13\ mM$. This decrease reflects the faster transition through the quiescent phase as glucose increases, consistent with the physiological response of beta cells, where higher glucose levels reduce the quiescent phase duration and increase bursting frequency. The linger time values ($35-5$ seconds) align with the typical range for  $\beta-$cells quiescent phases, as a slow burster, particularly at higher glucose levels. The longer linger time at $G = 7\ mM$ may indicate proximity to the pseudo-singular point, where a canard-induced delay prolongs the quiescent phase, a phenomenon that could be theoretically quantified using a weighted blow-up analysis.

\begin{figure}[ht]
    \centering
    \includegraphics[width=0.8\textwidth]{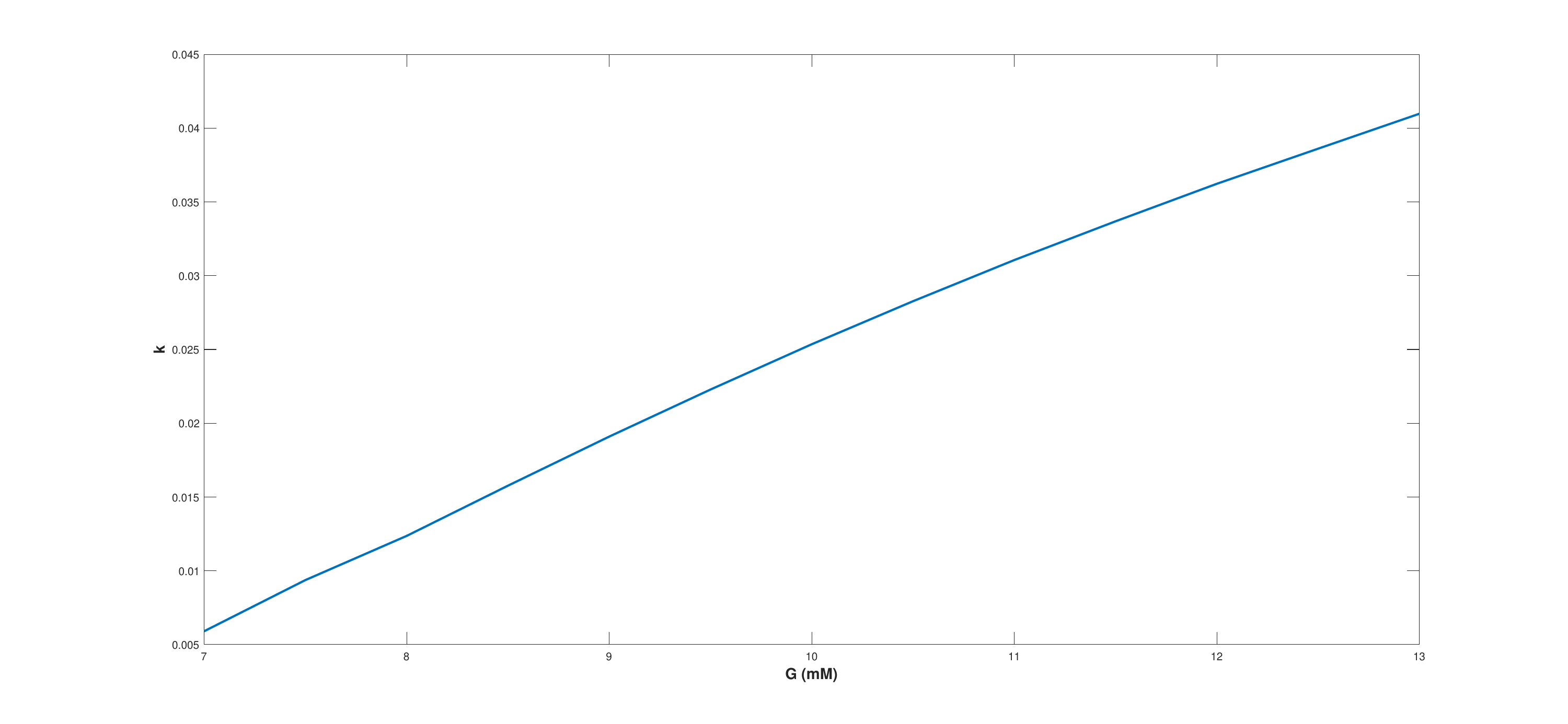} 
    \caption{Coupling strength ($k$) as a function of $t_{linger}$ for different glucose concentration ($G$).}
    \label{fig:k_vs_linger_time_G}
\end{figure}

Figure \ref{fig:k_vs_linger_time_G} represents the coupling strength $k$ as a function of glucose concentration ($ G $) for a local network of pancreatic $\beta-$cells, computed using the relation $k \propto \frac{1}{t_{linger}}$, with a desired error tolerance. The $x-$axis represents ($ G $) in millimolar ($mM$), ranging from $7$ to $13\ mM$. The $y-$axis shows $k$, ranging from $0.0028$ to $0.027$. The plot exhibits a smooth, monotonically increasing trend, with $k$ rising from $0.005$ at $G = 7\ mM$ to $0.042$ at $G = 13\ mM$. This increase corresponds to the decreasing linger time over the same ($G$) range, as stronger coupling is required to synchronize burst onsets within a shorter quiescent phase. The $k$ values are within the typical range for our studied $\beta-$cell model, ensuring synchronization within that phase. The trend highlights the glucose-dependent coupling needed for coordinated insulin release, with higher glucose levels necessitating stronger coupling to maintain synchronization amid faster bursting dynamics.

\section{Conclusion}\label{sec:conclusion}
Our study of pancreatic $\beta$-cell dynamics using a three-time-scale model has enlightened the mechanisms driving insulin secretion and their significance for glucose regulation and diabetes research. By examining the interplay of fast, intermediate, and slow variables, we explained how ATP-driven shifts in membrane potential produce the bursting patterns critical for insulin release. Through blow-up analysis, we identified canard-mediated transitions near pseudo-singular points that govern quiescent phase durations, decreasing from $35$ seconds at $7\ mM$ glucose to $5$ seconds at $13\ mM$. Numerical simulations further revealed that higher glucose concentrations reduce the bursting period from $10.3$ to $5.5$ minutes, requiring stronger coupling strengths ($0.005$ to $0.042$ ms$^{-1}$) to synchronize burst initiation across $\beta$-cells. These results underscore the pivotal role of synchronized oscillations in enhancing insulin pulses, providing insights into the impaired coordination observed in diabetes. By integrating mathematical rigor with biological insights, our work establishes a foundation for future studies on $\beta$-cell network dynamics, diabetes research, and potential therapeutic strategies for glucose homeostasis disorders.

\section*{Acknowledgments}
This work forms part of the ongoing PhD thesis ``Canard Dynamics and Synchronization in the Network of Three-Time-Scale Systems'' at Tohoku University. The author gratefully acknowledges Professor Hayato Chiba (AIMR, Tohoku University, Japan) for his invaluable guidance, continuous support, and patience throughout this research. This work was supported by JST University Fellowships for Science and Technology Innovation (Grant No. JPMJFS2102) and JST SPRING (Grant No. JPMJSP2114).

%\section*{Bibliography}

\end{document}